\documentclass[11pt, reqno]{amsart}
\pdfoutput=1

\usepackage{amssymb,amsmath,latexsym,amsthm, amsfonts}
\usepackage{epsfig}
\usepackage{mathrsfs}
\usepackage{graphicx}
\usepackage{epic,eepic,wrapfig,color,ifthen}
\usepackage{url}
\usepackage{hyperref}
\usepackage{pmboxdraw}
\usepackage{verbatim}
\usepackage{mathtools, textcomp}
\usepackage{cite}
\usepackage[normalem]{ulem}
\usepackage{enumerate}
\usepackage[bottom]{footmisc}
\usepackage{upgreek}

\newcommand{\cal}[1]{\mathcal{#1}}

\def\up{\upsilon}

\def\1{{\bf 1}}

\def\nn{\nonumber}

\def\d{{\mathsf{d}}}
\def\a{{\mathbf a}}
\def\b{{\mathbf b}}

 \def\sB {{\cal B}} \def\sC {{\cal C}}
\def\sD {{\cal D}}
 
\def\sF {{\cal F}}

\def\sL {{\cal L}}

 \def\sO {{\cal O}}

 \def\sU {{\cal U}}

\def\la{{\langle}}
\def\ra{{\rangle}}

\def \VDo {${\mathrm{VRD}}_{R_0}(U)$}

\def \VDi {$\mathrm{VRD}^{R_\infty}(\up)$}

\numberwithin{equation}{section}
\def\qed{{\hfill $\Box$ \bigskip}}

\def\FF{{\mathcal F}}

\def\GG{{\mathcal G}}

\def\R{{\mathbb R}}
\def\P{{\mathbb P}}
\def\N{{\mathbb N}}
\def\Z {{\mathbb Z}}
\def\E {{\mathbb E}} 
\def\LL {{\mathbb L}}
\def\X {{\mathfrak X}}

\def\vp{{\varphi}}
	
\def\tO {{\tilde{\Omega}}}
\def\tP {{\tilde{\P}}}

\def\fM {{\mathfrak B}_{b,+}(M)}
\def\fMg {{\mathfrak B}^\gamma_{b,+}(M)}
\def\fMi {{\mathfrak B}^\infty_{b,+}(M)}

\def\oB{{\overline B}}

\def\eps{\varepsilon}

\def\wt{\widetilde}
\def\pf{\noindent{\bf Proof.} }

\theoremstyle{plain}
\newtheorem{thm}{Theorem}[section]
\newtheorem{lem}[thm]{Lemma}
\newtheorem{cor}[thm]{Corollary}
\newtheorem{remark}[thm]{Remark}
\newtheorem{prop}[thm]{Proposition}
\newtheorem{defn}[thm]{Definition}

\newtheorem{example}[thm]{Example}

\theoremstyle{definition}
\newtheorem*{eg*}{Example}
\newtheorem*{egs*}{Examples}
\newtheorem*{def*}{Definition}

\theoremstyle{remark}

\begin{document}
\title[Laws of the iterated logarithm for occupation times of Markov processes]
{Laws of the iterated logarithm for occupation times of Markov processes}

\author{Soobin Cho, Panki Kim and  Jaehun Lee}\thanks{This work was supported by the National Research Foundation of Korea(NRF) grant funded by the Korea government(MSIP) 
(No. NRF-2021R1A4A1027378).
}

	\address[Cho]{Department of Mathematical Sciences,
		Seoul National University,
		Seoul 08826, Republic of Korea}
	\curraddr{}
	\email{soobin15@snu.ac.kr}
	
	\address[Kim]{Department of Mathematical Sciences and Research Institute of Mathematics,
		Seoul National University,
		Seoul 08826, Republic of Korea}

	\curraddr{}
	\email{pkim@snu.ac.kr}
	
	\address[Lee]{
	Korea Institute for Advanced Study,  
	Seoul 02455,
	Republic of Korea}

	\curraddr{}
	\email{hun618@kias.re.kr}

\maketitle

\begin{abstract}
	In this paper, we discuss the laws of the iterated logarithm (LIL) for occupation times of Markov processes $Y$ in general metric measure space
 both near zero and near infinity  under some minimal assumptions. 
We 
first establish LILs of (truncated)  occupation times  on balls $B(x,r)$ of radii $r$  up to an function $\Phi (r)$, which is an iterated logarithm of mean exit time  of $Y$, 
by showing that the function $\Phi$ is optimal.
Our first result on  LILs of occupation times covers
both near zero and near infinity regardless of transience and recurrence of the process. 
Our assumptions are truly local in particular at zero and the function $\Phi$ in 
our truncated occupation times $r \mapsto\int_0^{ \Phi (x,r)} \1_{B(x,r)}(Y_s)ds$ depends on space variable $x$ too.
We also prove that a similar LIL for total occupation times $r \mapsto\int_0^\infty \1_{B(x,r)}(Y_s)ds$ holds when the process is transient. Then we establish  LIL concerning large time behaviors of occupation times
$t \mapsto \int_0^t \1_{A}(Y_s)ds$
 under an additional condition that guarantees the recurrence of the process.
Our results cover a large class of Feller (Levy-like) processes, random conductance models with long range jumps,   jump processes with mixed polynomial local growths and jump processes with singular jumping kernels.

\medskip

\end{abstract}
\allowdisplaybreaks

\section{Introduction}
Law of the iterated logarithm (LIL), together with 
 the  the central limit theorem and law of large numbers,
 is considered as the fundamental limit theorem in Probability theory and a lot of beautiful results have been established for LIL for the sample path of Markov processes. See \cite{Bi86} and the references therein. 
 Very recently, in 
 \cite{CKL, CKL2}, the authors have obtained limsup and liminf laws for sample paths of a large class of standard Markov processes both near zero and near infinity under some minimal localized assumptions.

 On the other hand, 
the literature on LIL for occupation times for stochastic processes  is quite scarce and, to the best of our knowledge,  LIL for occupation times has been studied only for Brownian motion  and stable processes partially. See \eqref{e:intro_1}--\eqref{e:intro_3} below.

The purpose of this paper is to 
establish various forms of limsup LILs for occupation times of Markov processes in general metric measure space
  both near zero and infinity. 
  This paper is a continuation of our journey  on investigating minimal assumptions in the study of sample path properties for Markov processes via localizing approach.
  
  We first recall known results  on the limsup LIL for occupation times. Let $(Y_t)_{t \ge 0}$ be a Brownian motion  in $\R^d$, $d \ge 3$ and 
$B(0,r) \subset \R^d$ be the ball of radius $r$ centered at the origin.  
Denote by $\sU(B(0,r),t):=\int_0^t \1_{B(0,r)}(Y_s)ds$ the occupation time
on $B(0, r)$
 up to time $t$. Ciesielski and Taylor \cite{CT} proved that 
\begin{equation}\label{e:intro_1}
	\limsup_{r \to 0} \frac{\sU(B(0,r),\infty)}{r^2 \log |\log r|} = \frac{2}{p_d^2}, \quad \text{a.s.},
\end{equation}
where $p_d$ is the first positive zero of Bessel function $J_{d/2-2}(z)$ of the first kind. Using  the theory of large deviations,  Donsker and Varadhan \cite{DV2} showed a large scale counterpart of \eqref{e:intro_1}, namely, there exists a constant $c_1>0$ such that
\begin{equation}\label{e:intro_2}
	\limsup_{r \to \infty} \frac{\sU(B(0,r),\infty)}{r^2 \log \log r} = c_1,\quad \text{a.s.}.
\end{equation}
Motivated by \eqref{e:intro_2}, Shieh \cite{Sh} showed that for every strictly stable process $(Y_t)_{t \ge 0}$  in $\R^d$, $d\ge1$ with index $\alpha \in (0,2]$ (that is,    $Y_t$ and $s^{-1/\alpha}Y_{st}$ have the same distribution for all $s,t >0$) which has an everywhere strictly positive density, there exists a constant $c_2>0$ such that
\begin{equation}\label{e:intro_3}
	\limsup_{r \to \infty} \frac{\sU(B(0,r),r^2\log \log r)}{r^2 \log \log r} = c_2,\quad \text{a.s.}.
\end{equation}
Note that  only truncated occupation times considered in \cite{Sh} because $(Y_t)_{t \ge 0}$ can be recurrent in the setting of \eqref{e:intro_3} so that $\sU(B(0,r),\infty)=\infty$ almost surely.

In the first part of this paper, we study \eqref{e:intro_3}-type LIL and  its small scale  counterpart LIL for general standard (possibly recurrent) Markov processes   that have mere weak scaling structures. 
See Theorems \ref{t:generalocc-SP} and \ref{t:generalocc2} below. Then we get \eqref{e:intro_2}-type LIL under additional assumptions that guarantee the transience of the process. See Theorem \ref{t:generalocc3} below. Our assumptions on a scaling structure of processes are quite weak so that our result even covers some  random conductance model.

In the second part of the paper, we study large time behaviors of occupation times 
$t \mapsto \int_0^t \1_{A}(Y_s)ds$
and some additive functionals  when the process is recurrent. See Theorem \ref{t:occ1} below.

Our assumption is general enough to cover a large class of Levy-like processes,  jump processes with mixed polynomial local growths, jump processes with singular jumping kernels.
In particular, similar to 
 \cite{CKL, CKL2},
our assumptions at infinity allow  irregular behaviors of the process when it is on regions far away from the origin, which are controlled by a parameter $\up \in (0,1)$.
See Definition \ref{d:VD}(ii) and \textbf{Assumption B} below. Thanks to such weak assumptions, our results on LIL at infinity cover some random conductance model including  ones with long range jumps.  In Section \ref{s:Example}, we give two important examples in detail: (1) Feller processes on $\R^d$ and (2) Random conductance model.

\bigskip

\noindent {\bf Notations}: 
Without loss of generality, throughout this paper,  we use same fixed  positive constants $d_1$, $d_2$ and $C_i, i=0,1,2, \dots$ on conditions and statements both at zero and at infinity except in Section \ref{s:Example}. Lower case letters  $a_i$ and $c_i$, $i=0,1,2,...$, which denote positive real constants, are fixed in each statement and proof and  the labeling of these constants starts anew in each proof unless they are specified to denote particular values. 
We use the symbol ``$:=$'' to denote a definition, 
which is read as ``is defined to be.''  
We write $a\wedge b:=\min\{a,b\}$, $a\vee b:=\max\{a,b\}$ and $\lfloor a \rfloor := \sup\{ n \in \Z : n \le a \}$.
$\overline{A}$ denotes the closure of $A$. We extend a function $f$ defined on $M$ to $M_\partial$ by setting $f(\partial)=0$. The notation $f(x) \simeq g(x)$ means that there exist constants $c_2 \ge c_1>0$ such that $c_1g(x)\leq f (x)\leq c_2 g(x)$ for a specified range of the variable $x$.

\section{Settings and Main results}\label{s:setting}

\subsection{Settings}
Throughout this paper, we assume that $(M,d)$ is a locally compact separable metric space with a base point $o \in M$, and $\mu$ is a positive Radon measure on $M$ with full support. Denote by $B(x,r)$ the open ball centered at $x$ with radius $r$, and $\sB(M)$ the family of all Borel sets on $M$.  

For  $D \in \sB(M)$, we write
\begin{equation*}
	\updelta_D(x)=\inf \{ d(x,y) : y \in M \setminus D\}.
\end{equation*}
Define
\begin{equation}\label{e:def_d}
	\d (x) = d(x,o) + 1.
\end{equation}
Note that when $(M,o)=(\R^d, 0)$,  $\d(x)$ is equal to $|x| + 1$.
Since $\d(x) \ge 1$,  the map $\up \to \d(x)^\up$ is non-decreasing on $(0,\infty)$ for all $x \in M$. The function $\updelta_D$ will be used to describe assumptions and statements for LILs near zero, and $\d$ will be used for LILs near infinity.

Write  $V(x,r)=\mu(B(x,r))$. We recall interior and  local and weak versions of \textit{the volume doubling and reverse doubling property}  from the authors' previous paper \cite{CKL}.   

\begin{defn}\label{d:VD}
	{\rm (i) For an open set $U \subset M$ and $R_0 \in (0,\infty]$, we say that 
	\textit{an interior volume doubling and reverse doubling property} 
	\VDo \  holds if there exist  constants $C_0 \in (0,1)$ and $d_1,d_2,c_\mu,C_\mu>0$ such that for all $x \in U$ and $0<s\le r<R_0 \wedge (C_0 \updelta_U(x))$,
		\begin{equation}\label{e:VRD}
			c_\mu \left( \frac{r}{s} \right)^{d_1} \le \frac{V(x,r)}{V(x,s)} \le  C_\mu \left( \frac{r}{s} \right)^{d_2}.
		\end{equation}
		
		\noindent (ii) 	For  $R_\infty\ge1$ and $\up \in (0,1)$, we say that \textit{a weak volume doubling and reverse doubling property at infinity}  \VDi \ holds  if there exist constants $d_1,d_2, c_\mu, C_\mu>0$  such that  \eqref{e:VRD} holds for all $x \in M$ and $r \ge s > R_\infty \d(x)^\up$.  
} \end{defn}

Let $X = (\Omega, \sF_t, X_t,  t \ge 0;  \P^x, x \in M_\partial)$ be  a  Borel standard Markov process on $M_\partial:=M \cup \{\partial\}$ where $\partial$ is a cemetery point added to $M$.  We sometimes denote $X(t)$ for $X_t$. 
According to \cite[Theorem 1.1]{BJ73}, every Borel standard Markov process on $M_\partial$ has a L\'evy system. In  this paper, we   assume that $X$ has a   L\'evy system of the form $(J(x,\cdot),ds)$, that is, there exists a family of Borel measures  $(J(x,dy):x \in M)$ on $M_\partial$ such that   for any  $z \in M$, $t>0$ and any non-negative Borel function $F$ on $M \times M_\partial$ vanishing on the diagonal,
\begin{equation*}
	\E^z \bigg[ \,\sum_{s \le t} F(X_{s-}, X_s)\,  \bigg] = \E^z \bigg[ \int_0^t \int_{M_\partial} F(X_s,y)J(X_s,dy)ds \bigg].
\end{equation*}
The measure $J$  is called the \textit{L\'evy measure} of  $X$.   Note that we count the killing term $J(x, \partial)$ as a part of the L\'evy measure.

For $D\in \sB(M)$, we denote  \textit{the first exit time} of $X$  from $D$ by
$$
\tau_D=\tau[D]:= \inf \{ t>0 : X_t \in M_\partial \setminus D \}.
$$
It is known that the part process $X^D$ of $X$ in $D$, which is defined by   $X_t^D := X_t \, \1_{\{\tau_D > t\}} + \partial \, \1_{\{\tau_D \le t\}}$,  is a Borel standard process having a L\'evy system induced from that of $X$. See \cite[Section 3.3]{CF} for details.  We call  a Borel  measurable  function $p^D: (0,\infty) \times D \times D \to [0,\infty]$ is the  \textit{heat kernel} of   $X^D$ if it satisfies the following properties:

\smallskip

\setlength{\leftskip}{4mm}
\noindent (1) $\E^x [f(X^D_t)] = \int_{D} p^D(t,x,y)f(y) \mu(dy)$ for all  $f \in L^\infty(D;\mu)$, $t>0$ and $x \in D$;

\smallskip

\noindent (2) $p^D(t+s,x,y) = \int_{D} p^D(t,x,z) p^D(s,z,y) \mu(dz)$ for all $t,s>0$ and $x,y \in D$.

\smallskip

\setlength{\leftskip}{0mm}

\noindent 
We write  $p(t,x,y)$ for $p^M(t,x,y)$.

\medskip

Now, we introduce  our assumptions; cf. \cite{CKL2}. Fix an open subset $U$ of $M$.

\medskip

\noindent \textbf{Assumption A.} There exist  constants $R_0>0$, $C_0\in(0,1)$, $C_1>1$, $C_i>0$, $2\le i\le 7$  and a function $\phi:U \times (0,R_0) \to (0,\infty)$ such that $\phi(z,\cdot)$ is increasing and continuous for every $z \in U$ and the following properties are satisfied  for every $x \in U$ and  $0<r<R_0 \wedge (C_0\updelta_U(x))$:
\begin{equation}\label{A1}\tag{${\rm CE}_{R_0}(\phi,U)$}
	C_1^{-1} \phi(y,r)\le  \phi(x,r) \le C_1 \phi(y,r)\quad \text{for all} \; y \in B(x,r);
\end{equation}\\[-7mm]
\begin{equation}\label{A2}\tag{${\rm US}_{R_0}(\phi,U)$}
	\lim_{r \to 0} \phi(x,r) = 0, \quad  \phi(x,r) \le C_2 \phi(x,r/2);
\end{equation}\\[-7mm]
\begin{equation}\label{A3}\tag{${\rm Tail}_{R_0}(\phi, U,\le)$}
	J(x,M_\partial \setminus B(x,r)) \le  \frac{C_3}{ \phi(x,r)};
\end{equation}\\[-5mm]
\begin{equation}\label{A4}\tag{${\rm SP}_{R_0}(\phi,U)$}
C_4e^{-C_5 n} \le \P^x\big(\tau_{B(x,r)} \ge n \phi(x,r)\big)\ge 	C_5e^{-C_7 n}  \quad \text{for all} \; n \ge 1.
\end{equation}

\bigskip
The most delicate part of Assumption A is \ref{A4}.
We will see in Proposition \ref{p:E} that under \VDo \ and \ref{A2}, the following condition implies \ref{A4}:  There exist constants $c_l>0$ and $\eta  \in (0,1)$ such that for every $x \in U$ and $0<r<R_0 \wedge (C_0 \updelta_U(x))$, the heat kernel $p^{B(x,r)}(t,x,y)$ of $X^{B(x,r)}$ exists and satisfies that
\begin{equation}\label{A4+}\tag{${\rm NDL}_{R_0}(\phi,U)$}
\hspace{0.3in}	p^{B(x,r)}( \phi(x, \eta r),y,z) \ge \frac{c_l}{V(x, r)} \quad \text{for all} \;\;  y,z \in B(x, \eta^2r).
\end{equation}

\medskip

Next, we give assumptions for LILs at infinity. 

\medskip

\noindent \textbf{Assumption B.} There exist constants $R_\infty\ge 1$, $\up,\eta\in(0,1)$, $c_l>0$,  $\ell>1$,  $C_i>0$, $2\le i\le 3$ and an increasing continuous function $\varPhi:(R_\infty,\infty) \to (0,\infty)$ such that the following properties are satisfied  for every $x \in M$ and  $r>R_\infty \d(x)^\up$:
\begin{equation}\label{B2}\tag{${\rm ULS}^{R_\infty}(\varPhi)$}
	2 \varPhi(s/\ell) \le  \varPhi(s) \le C_2	\varPhi(s/2)  \quad \text{for all} \; s >R_\infty,
\end{equation}\\[-7mm]
\begin{equation}\label{B3}\tag{${\rm Tail}^{R_\infty}(\varPhi,\up,\le)$}
	J(x,M_\partial \setminus B(x,r)) \le  \frac{C_3}{\varPhi(r)},
\end{equation}\\[-7mm]
\begin{align}
&\textit{the heat kernel $p^{B(x,r)}(t,x,y)$ of $X^{B(x,r)}$ exists and satisfies that}\nn\\
&	p^{B(x,r)}( \varPhi(\eta r),y,z) \ge \frac{c_l}{V(x, r)} \quad \text{for all} \;\;  y,z \in B(x, \eta^2r).\label{B4+}\tag{${\rm NDL}^{R_\infty}(\varPhi,\up)$} 
\end{align}

\medskip

\begin{remark}\label{r:ass-B}
	{\rm 	(i) The upper scaling condition at zero \ref{A2} is equivalent with that there exist constants $\beta_2>0$ and $C_U\ge 1$ such that
		\begin{equation}\label{e:scaling-zero}
			\frac{\phi(x,r)}{\phi(x,s)} \le C_U \bigg(\frac{r}{s}\bigg)^{\beta_2} \quad \text{for all} \; x \in U, \; 0<s\le r<R_0 \wedge (C_0\updelta_U(x)),
		\end{equation}
		and the upper and lower scaling condition at infinity \ref{B2} is equivalent with that  there exist constants $\beta_1,\beta_2>0$ and $C_U\ge 1 \ge C_L>0$ such that
		\begin{equation}\label{e:scaling-infty}
			C_L \bigg(\frac{r}{s}\bigg)^{\beta_1} \le 	\frac{\varPhi(r)}{\varPhi(s)} \le C_U \bigg(\frac{r}{s}\bigg)^{\beta_2} \quad \text{for all} \;  r \ge s> R_\infty.
	\end{equation}
In particular,  \ref{B2} implies $\lim_{r \to \infty} \varPhi(r)=\infty$. Note that we do not impose lower scaling condition at zero.

\noindent (ii) One can check that under \ref{A4}, there exists a constant $\Lambda \ge 1$ such that
\begin{equation}\label{e:mean-exit}
	\Lambda^{-1} \E^x[\tau_{B(x,r)}] \le \phi(x,r) \le 	\Lambda \E^x[\tau_{B(x,r)}]
\end{equation}
 for every $x \in U$ and  $0<r<R_0 \wedge (C_0\updelta_U(x))$.
 See  \cite[(4.12)--(4.13)]{CKL}. 
We will see that  \VDi, \ref{B2} and  \ref{B4+} imply a near infinity counterpart of \ref{A4}. Hence,  under these assumptions, there exist  constants $R_\infty',\Lambda\ge 1$ such that \eqref{e:mean-exit} holds with $\varPhi(r)$ instead of $\phi(x,r)$ for every $x \in M$ and $r>R_\infty'\d(x)^\up$. See Proposition \ref{p:E}(ii).
}
\end{remark}

We recall the following remark from the previous paper \cite{CKL2}.

\begin{remark}\label{r:conditions}
	{\rm (i) We impose conditions at infinity  only for $r>R_\infty \d(x)^\up$. By considering such weak assumptions at infinity, our results cover  some random conductance models. See Subsection \ref{s:RCM}.

		\noindent (ii) We will use a proper zero-one law for tail events (Proposition \ref{p:law01}) to get a deterministic LILs at infinity. To use this zero-one law, we need \ref{B4+}, which is stronger than a near infinity counterpart of \ref{A4}, for LILs at infinity.
	}
\end{remark}

\subsection{Main results}

For a Borel measurable function $F$ on $M$ and $t \in [0,\infty]$, we denote $\sU(F,t)$ for the \textit{(truncated) occupation time} of $X$ for $F$ until $t$, namely,
\begin{equation*}
	\sU(F, t)(\omega):= \int_0^t F(X_s(\omega))ds.
\end{equation*}
We write $\sU(D,t):=\sU(\1_D,t)$ for $D \in \sB(M)$.

The first set of LILs concerns occupation times in balls. 

\begin{thm}\label{t:generalocc-SP}
	Suppose that \ref{A1}, \ref{A2}, \ref{A3} and \ref{A4} hold for an open set $U \subset M$. Let $f_0:U \times (0,\infty) \to [0,1]$ be a deterministic function such that for every $x \in U$ and $\kappa>0$, 
	\begin{equation}\label{e:occ_SP}
		\limsup_{r \to 0} \frac{\sU\big(B(x, r),\, \kappa \phi(x,r) \log |\log \phi(x,r)|\big)}{\kappa \phi(x,r) \log |\log \phi(x,r)|} = f_0(x, \kappa) ,~~\quad
		\P^x\mbox{-a.s.},
	\end{equation}
which is well-defined by the Blumenthal's zero-one law. 	Then $f_0$ satisfies the following properties:
	
	{\rm (P1)} For each fixed $x \in U$, the map $\kappa \mapsto f_0(x,\kappa)$ is strictly positive and non-increasing.
	
		{\rm (P2)} There exists a constant $\kappa_0>0$ such that $f_0(x,\kappa) = 1$ for all $x \in U$ and $\kappa \le \kappa_0$.
	
		{\rm (P3)} $\lim_{\kappa \to \infty} f_0(x,\kappa) = 0$ uniformly on $U$.
\end{thm}

Properties (P1)-(P3) of $f_0$ in the above theorem and \eqref{e:finfty} below show that the results in Theorems \ref{t:generalocc-SP} and \ref{t:generalocc2} are  optimal. 
See Remark \ref{r:optimal} below.
The proofs of (P3) in Theorems \ref{t:generalocc-SP}  and \eqref{e:finfty} are the most delicate parts of the paper.  See Lemmas  \ref{l:generalocc-SP-4}--\ref{l:generalocc-SP-3} below.

\medskip

By Proposition \ref{p:E}(i) below,  Theorem  \ref{t:generalocc-SP} implies 

\begin{cor}\label{c:generalocc}
Suppose that \VDo, \ref{A1}, \ref{A2}, \ref{A3} and \ref{A4+} hold for an open set $U \subset M$. Then the function $f_0:U \times (0,\infty) \to [0,1]$ defined by  \eqref{e:occ_SP}  satisfies  {\rm (P1)-(P3)} in Theorem \ref{t:generalocc-SP}.  
\end{cor}

\medskip

To obtain the (uniform) limit of occupation time at infinity, we introduce a regularized function
\begin{equation}\label{e:molli1}
\vp(r):= \int_1^r s^{-1}\varPhi(s)ds, \qquad r>1.
\end{equation}
See \eqref{e:compareocc} below for a reason of this regularizing. Under \ref{B2}, there is  a constant $C>1$ such that  
\begin{equation}\label{e:molli2}
C^{-1}\varPhi(r) \le \vp(r) \le C \varPhi(r) \quad \text{ for all} \;\, r \ge 2.
\end{equation}
Indeed, by the monotone property of $\varPhi$ and  \eqref{e:scaling-infty}, it holds that for any  $r \ge 2$,  
$$2^{-\beta_2}C_U^{-1} \log 2\le \frac{\varPhi(r/2)}{\varPhi(r)} \int_{r/2}^r s^{-1}ds \le \frac{\vp(r)}{\varPhi(r)} \le  C_L^{-1}r^{-\beta_1} \int_1^r s^{-1+\beta_1}ds  \le \beta_1^{-1}C_L^{-1}.$$
By  \eqref{e:scaling-infty} and \eqref{e:molli2}, there exist constants $C_U' \ge 1  \ge C_L'>0$ such that 
\begin{equation}\label{e:scaling-infty*}
 C_L' \bigg(\frac{r}{s}\bigg)^{\beta_1}\le	\frac{\vp(r)}{\vp(s)} \le C_U' \bigg(\frac{r}{s}\bigg)^{\beta_2} \quad \text{for all} \;  r \ge s\ge 2R_\infty.
\end{equation}

\begin{thm}\label{t:generalocc2}
	  Suppose that \VDi, \ref{B2}, \ref{B3} and \ref{B4+} hold.  Then  there exists 
	  a deterministic non-increasing strictly positive function $f_\infty:(0,\infty) \to (0,1]$  such that for every $x \in M$ and  $\kappa>0$,
	\begin{equation}\label{e:generalocc2_1}
	\limsup_{r \to \infty} \frac{\sU\big(B(x, r),\, \kappa\vp(r) \log \log \vp(r)\big)}{\kappa \vp(r) \log \log \vp(r)} = f_\infty(\kappa),~~\quad
	\P^z\mbox{-a.s.}, ~~~\forall z \in M,
	\end{equation}
where  $\vp$ is defined by \eqref{e:molli1}. 
Moreover, there exists a constant $\kappa_\infty>0$ such that 
\begin{align}
\label{e:finfty}
f_\infty(\kappa) = 1 \;\;\text{for }\kappa \le \kappa_\infty \quad  \text{and} \quad \lim_{\kappa \to \infty} f_\infty(\kappa)=0.
\end{align}
\end{thm}

\begin{remark}\label{r:optimal}
	{\rm  Let $\psi :U \times (0,\infty) \to (0,\infty)$ be an increasing function such that for some $x_0 \in U$,
	$$ \lim_{r \to 0} \frac{\psi(x_0,r)}{\phi(x_0,r) \log |\log \phi(x_0,r)|} = 0 \quad  (resp. \,=\infty). $$
	Then under the setting of Theorem \ref{t:generalocc-SP}, by (P1)-(P3) and \eqref{e:occ_SP}, it holds that    for any $\kappa>0$,
	\begin{align*}
	\limsup_{r \to 0} \frac{\sU(B(x_0,r),\kappa \psi(x_0,r))}{\kappa \psi(x_0,r)} &=  \limsup_{r \to 0} \frac{\sU\Big(B(x_0,r), \big( \kappa  \frac{\psi(x_0,r)}{\phi(x_0,r) \log |\log \phi(x_0,r)|} \big) \phi(x_0,r) \log |\log \phi(x_0,r)|\Big)}{\big( \kappa  \frac{\psi(x_0, r)}{\phi(x_0,r) \log |\log \phi(x_0,r)|} \big) \phi(x_0,r) \log |\log \phi(x_0,r)|} \\
	&= \limsup_{r \to 0}\, f_0\big(x_0,\kappa  \frac{\psi(x_0,r)}{\phi(x_0,r) \log |\log \phi(x_0,r)|}\big) = 1 \quad (resp. \,=0).
	\end{align*}
Thus, $\psi$ can not satisfy both (P2) and (P3) simultaneously. 
Analogous result holds concerning $\vp(r)$ instead of $\phi(x_0,r)$.
 In this sense, the rate functions $\phi(x,r)\log|\log \phi(x,r)|$ and $\vp(r)\log\log\vp(r)$ in \eqref{e:occ_SP} and \eqref{e:generalocc2_1} are optimal and unique.}
\end{remark}

\medskip

Next, we consider limsup behaviors of total occupation times $\sU(B(x,r),\infty)$ when the process $X$ is transient.  
To state this result, we need the following near diagonal upper heat kernel estimates.
\begin{align}
	&\hspace{-3mm}\textit{The heat kernel $p(t,x,y)$ of $X$ exists and there exists a constant $c_u>0$ such that}\nn \\
	&\hspace{-3mm}p(t,x,y) \le  \frac{c_u}{V(x,\varPhi^{-1}(t)) \wedge V(y,\varPhi^{-1}(t))} \;\;\text{ for all} \;\, x,y \in M, \, t\ge  \varPhi\big(R_\infty (\d(x) \vee \d(y))^\up\big). \label{NDU}\tag{${\rm NDU}^{R_\infty}(\varPhi,\up)$}
\end{align}

\begin{remark}
	{\rm Suppose that $X$ is $\mu$-symmetric and has the heat kernel $p(t,x,y)$.  If there exist $R_\infty \ge1$, $\up \in (0,1)$ and  $c_1>0$  such that 
		\begin{equation}\label{e:point_wDUi}
			p(t,x,x) \le \frac{c_1}{V(x,\varPhi^{-1}(t))} \quad\text{for all} \;\; x \in M, \,  t> \varPhi(R_\infty\d(x)^\up),
		\end{equation}
		then \ref{NDU} is satisfied. Indeed, by the semigroup property, $\mu$-symmetry  and  Cauchy–Schwarz inequality, \eqref{e:point_wDUi} implies that  for all $x,y \in M$ and $t \ge \varPhi(R_\infty (\d(x) \vee \d(y))^\up)$, 
		\begin{align*}
			p(t,x,y)&= \int_M p(t/2,x,z)p(t/2,z,y)\mu(dz) \le \left(\int_M p(t/2,x,z)^2\mu(dz) \right)^{1/2}  \left(\int_M p(t/2,y,z)^2\mu(dz) \right)^{1/2}  \\
			& = p(t,x,x)^{1/2}p(t,y,y)^{1/2}\le p(t,x,x) \vee p(t,y,y) \le \frac{c_1}{V(x,\varPhi^{-1}(t))\wedge V(y,\varPhi^{-1}(t))}.
		\end{align*}
	}
\end{remark}

\begin{thm}\label{t:generalocc3}
Suppose that \VDi, \ref{B2}, \ref{B3}, \ref{B4+} and \ref{NDU} hold. Suppose also that the upper inequality in \eqref{e:scaling-infty} holds with  $\beta_2<d_1$, where $d_1$ is the constant from \VDi. Let $f_\infty$ be the function given in Theorem \ref{t:generalocc2}. Then the limit $\kappa_1:=\lim_{\kappa \to \infty} \kappa f_\infty(\kappa)$ exists in $(0,\infty)$ and for every $x \in M$,
	\begin{equation}\label{e:occ3}
		\limsup_{r \to \infty} \frac{\sU\big(B(x, r),\,\infty\big)}{ \vp (r) \log \log \vp(r)} =\kappa_1,~~\quad
		\P^z\mbox{-a.s.}, \;\;\forall z \in M.
	\end{equation}
\end{thm}

\medskip

As we mentioned in the introduction, using the theory of large deviations, Donsker and Varadhan \cite{DV2} showed \eqref{e:occ3} for transient Brownian motions, and Shieh \cite{Sh} showed \eqref{e:generalocc2_1} for every strictly stable process $X_t$ in $\R^d$ $(d\ge1)$ that has an everywhere strictly positive density. But their proofs use the  strict scaling property of the process which is  not available in the present paper.
In this paper,  we take a different approach relying on our Theorem \ref{t:generalocc2} and a Chung-type liminf LIL.

Lastly, we concern large time behaviors of $t \mapsto \sU(F,t)$ when $X$ is recurrent (see Remark \ref{r:NDLrecur}). We define a function $\Theta$ on $M \times [0, \infty) \times [0,\infty]$ as
\begin{equation}\label{deftheta}
	\Theta(x,r,t)=\Theta(t,r,t;\vp):=\int_{\vp(r)}^t \frac{ds}{V(x, \vp^{-1}(s))}.
\end{equation}
We simply denote $\Theta(t)$ for $\Theta(o,1,t)$. See Lemma \ref{l:theta} for basic properties of $\Theta$. 

Denote by $\fM$ the family of all bounded non-negative Borel measurable functions on $M$.  Recall the definition of $\d$ from \eqref{e:def_d}.  For $\gamma \in (0,\infty]$, we set
\begin{equation}\label{e:Linfty}
\fMg:= \begin{cases}
\big\{F \in  \fM: \d(x)^\gamma F(x) \in \fM \big\} \quad &\mbox{if} \;\; \gamma<\infty;\\[3pt]
\big\{F \in  \fM: F \text{ is compactly supported}  \big\} \quad &\mbox{if} \;\; \gamma=\infty.
\end{cases}
\end{equation}
Under \VDi, we have  $\fMg \subset L^1(M,\mu) \cap L^\infty(M,\mu)$ if $\gamma \in (d_2,\infty]$ where $d_2$ is the constant from \eqref{e:VRD}. Indeed, for every $F \in \fMg$ with $\gamma \in (d_2,\infty]$, we have
\begin{align*}
	&\int_M F(y) \mu(dy) \le   \sum_{n=0}^\infty \int_{\{y \in M: 2^n \le \d(y) < 2^{n+1}\}} 2^{-n\gamma }\d(y)^{\gamma }F(y) \mu(dy)\\
	&\le c_1 
	\lVert \d^\gamma F
	 \rVert_{L^\infty(M,\mu)} \sum_{n=0}^\infty  2^{-n\gamma } V(o, 2^{n+1})  \le c_2
	  \lVert	  \d^\gamma F
	   \rVert_{L^\infty(M,\mu)} V(o,1) \sum_{n=1}^\infty  2^{-n(\gamma -d_2)} <\infty.
\end{align*} 

We suppose that one of the following two conditions holds true:
\begin{equation}\label{e:F}
(1) \; F \in \fMi \quad \;\text{ or }\; \quad
(2) \; \text{$\beta_1 > d_2$ and $F \in \fMg$ for some $\gamma \in \displaystyle\big( \frac{\beta_1 d_2}{\beta_1 - d_2} , \infty \big).$}
\end{equation}
Note that each of \eqref{e:F} implies  $\lVert F \rVert_{L^1(M,\mu)}<\infty$.

\begin{thm}\label{t:occ1}
	Suppose that  \VDi,  \ref{B2}, \ref{B3}, \ref{B4+} and \ref{NDU} hold.  Suppose also that
		\begin{equation}\label{e:recur}
	 \Theta(x_0, 1,  \infty) = \infty \quad \text{for some} \;\, x_0 \in M.
	\end{equation}
Then  there are constants $0<a_1 \le a_2<\infty$ such that for every $F \in \fM$ satisfying \eqref{e:F} and $\Vert F \Vert_{L^1(M,\mu)} \neq 0$, there exists a constant $a_F \in [a_1, a_2]$ satisfying
	\begin{equation}\label{e:occ_generalform}
	\limsup_{t \to \infty} \frac{\sU(F,t)/ \Vert F \Vert_{L^1(M,\mu)}}{\Theta\big(t/\log \log \Theta(t)\big) \log \log  \Theta(t)} = a_F,\quad\;\;
	\P^z\mbox{-a.s.},~~\forall z\in M.
	\end{equation}
\end{thm}
\medskip
When the index  in the lower scaling condition on $\varPhi$ at infinity is bigger than 
the upper dimension in \VDi, the function $\Theta$ can be written explicitly without integral form so that we get the following corollary.

\begin{cor}\label{c:occ1}
	Suppose that  \VDi,  \ref{B2}, \ref{B3}, \ref{B4+} and \ref{NDU} hold. Suppose also that the lower inequality in \eqref{e:scaling-infty} holds with  $\beta_1>d_2$, where $d_2$ is the constant from \VDi. Then \eqref{e:recur} is satisfied and  there are constants $0<\wt a_1 \le \wt a_2<\infty$ such that  for every  $F \in \fM$ satisfying \eqref{e:F} and $\Vert F \Vert_{L^1(M,\mu)} \neq 0$, there exists a constant $\wt a_F \in [\wt a_1, \wt a_2]$ satisfying 
	\begin{equation}\label{e:occ_simpleform}
		\limsup_{t \to \infty} \frac{\sU(F,t)/\Vert F \Vert_{L^1(M,\mu)}}{ t/V\big(o, \vp^{-1}(t/\log \log t)\big)} = \wt a_F,\quad\;\;
		\P^z\mbox{-a.s.},~~\forall z\in M.
	\end{equation}
\end{cor}

\begin{remark}\label{r:NDLrecur}
{\rm Suppose that  \VDi, \ref{B2}, \ref{B4+} and \ref{NDU} hold. Then   \eqref{e:recur} is satisfied  if and only if  $X$ is \textit{Harris recurrent}, namely:
\begin{equation}\label{e:harris}
\text{For all $z \in M$ and $A \in \sB(M)$ with $\mu(A)>0$, it holds that $\P^z\big(\,\sU(A, \infty) = \infty \big) = 1$.}
\end{equation}
In particular, \eqref{e:recur} holds if and only if $\Theta(x,1,\infty)=\infty$ for all $x \in M$.

\medskip

\noindent \textbf{Proof of Remark \ref{r:NDLrecur}.} Suppose that \eqref{e:harris} holds.  Let $x_0 \in M$. By \cite[Proposition 2.4]{Ge80},  \eqref{e:harris} yields  that $\int^\infty p(t, x_0, y) dt = \infty$ for $\mu$-a.e. $y \in M$. For any  $y \in B(x_0, \d(x_0))$, it holds that
$R_\infty(\d(x_0) \vee \d(y))^\up  \le 2R_\infty \d(x_0)$ by the triangle inequality. Hence, using \ref{NDU}, \VDi \ and \eqref{e:molli2}, we get that for all  $t > \phi(2R_\infty \d(x_0))$ and $y \in B(x_0, \d(x_0))$, 
\begin{equation*}
	p(t,x_0, y) \le \frac{c_u}{V(x_0, \phi^{-1}(t)) \wedge V(y, \phi^{-1}(t))}\le \frac{c_u}{V(x_0, \phi^{-1}(t)/2)} \le \frac{c_1}{V(x_0, \vp^{-1}(t))}.
\end{equation*}
In the second inequality above, we used the fact that $B(x_0, \phi^{-1}(t)/2) \subset B(y, d(x,y) + \phi^{-1}(t)/2 ) \subset B(y, \phi^{-1}(t))$  for all  $t > \phi(2R_\infty \d(x_0))$ and $y \in B(x_0, \d(x_0))$. Therefore, 
by \eqref{e:molli2}, there is $c_2>1$ and  $y \in B(x_0,\d(x_0))$  such that
$$
\infty = \int_{c_2\phi(2R_\infty \d(x_0))}^\infty p(t,x_0, y) dt \le \int_{\vp(2R_\infty \d(x_0))}^\infty \frac{c_1}{V(x_0, \vp^{-1}(t))}dt \le  c_1 
\Theta(x_0, 1, \infty).$$

Conversely, suppose that  \eqref{e:recur} holds for some $x_0 \in M$. Let $y,z \in M$. By \ref{B4+}, \eqref{e:scaling-infty} and \eqref{e:molli2}, there exist  $T,c_1,c_2>0$ such that $p(t,y,z) \ge p^{B(x_0, c_1\phi^{-1}(t))}(t, y,z) \ge c_2/V(x_0,\phi^{-1}(t))$ for all $t>T$. Hence, $\int^\infty p(t,y,z)dt  = \infty$ by \eqref{e:molli2} and \eqref{e:recur}. Now we deduce from  \cite[Proposition 2.4]{Ge80} and \cite[Theorem 1]{KM} that \eqref{e:harris} holds true. \qed
}
\end{remark}

 \section{Preliminary}

 Recall that $\Theta$ is defined in  \eqref{deftheta} and that $\vp$ is defined in  \eqref{e:molli1}.
  We begin this section with some basic properties of  $\Theta$.
 
 \begin{lem}\label{l:theta}
 	(i) For every $x \in M$, it holds that
 	\begin{equation*}
 	\Theta(x,r,t;\vp)  \ge \frac{t}{2V(x,\vp^{-1}(t))}\quad\;\; \text{for all} \;\; t \ge 2\vp(r).
 	\end{equation*}

 \noindent (ii) For every $x \in M$, it holds that
 \begin{equation*}
 	\frac{\Theta(x,r,t)}{\Theta(x,r,u)} \le 3 \left(\frac{t}{u}\right)^{\log 3/ \log 2}  \quad \text{for all} \;\; t>u \ge 2\vp(r).
 \end{equation*}

 	\noindent (iii) Suppose that  \VDi \ holds. If \eqref{e:scaling-infty} holds with $\beta_1>d_2$,  where $d_2$ is the constant from \eqref{e:VRD}, then there exists $c_1>0$ such that for all $x \in M$, $r > R_\infty \d(x)^\up$ and $t \ge 2\vp(r)$,
 	\begin{align*} \frac{t}{2V(x,\vp^{-1}(t))}\le 	\Theta(x,r,t) \le  \frac{c_2t}{V(x,\vp^{-1}(t))}.
 	\end{align*}
 
 \end{lem}
 \pf  (i) Using the monotone property of $\vp$, we get that for all $t \ge 2\vp(r)$,
 $$\Theta(x,r,t) \ge \int_{t/2}^{t} \frac{ds}{V(x, \vp^{-1}(s))}\ge \frac{t}{2V(x, \vp^{-1}(t))}.$$

 \noindent (ii)  Fix $t>u \ge 2\vp(r)$ and let
 $n_0\ge 1$ be the smallest integer such that $2^{n_0} u \ge t$.   Using the monotone property of $\vp$ and (i), we get that 
 $$\Theta(x,r,2u) - \Theta(x,r,u) = \int_{u}^{2u} \frac{ds}{V(x, \vp^{-1}(s))} \le \frac{u}{V(x, \vp^{-1}(u))} \le 2 \Theta(x,r,u).$$
 By iterating the above inequalities, we deduce that
 \begin{equation*}
 \Theta(x,r,t) \le \Theta(x,r,2^{n_0}u) \le 3\Theta(x,r,2^{n_0-1}u) \le \cdots \le   3^{n_0} \Theta(x,r,u) \le 3 (t/u)^{\log 3/\log 2} \Theta(x,r,u).
 \end{equation*} 
 
 \noindent (iii) Using \VDi \ and \eqref{e:scaling-infty*}, we obtain that for all $x \in M$, $r>R_\infty \d(x)^\up$ and $t \ge 2\vp(r)$,
 \begin{align*}
 	\Theta(x,r,t) &\le \frac{c_1}{V(x,  \vp^{-1}(t))} \int_{ \vp(r)}^t \left( \frac{ \vp^{-1}(t)}{ \vp^{-1}(s)}\right)^{d_2} ds  \le \frac{c_2}{V(x,  \vp^{-1}(t))} \int_{ \vp(r)}^t \left( \frac{t}{s}\right)^{d_2/\beta_1} ds \\
 	& \le \frac{c_2}{V(x,  \vp^{-1}(t))} \int_{0}^t \left( \frac{t}{s}\right)^{d_2/\beta_1} ds = \frac{c_2 \beta_1 }{(\beta_1-d_2)} \frac{1 }{V(x,  \vp^{-1}(t))}.
 \end{align*}
\qed

For $A \in \sB(M)$, denote by $\sigma_A= \inf\{t > 0 : X_t \in A\}$ the first hitting time of $A$.

\begin{lem}\label{l:suprema}
Let $F$ be a non-negative Borel measurable function on $M$ and $t \in [0,\infty]$. Suppose that supp$[F] \subset A$ for some close set $A$. Then it holds that
$$
\sup_{x \in A} \E^x\big[\,\sU(F,t)\big] = \sup_{x \in M} \E^x\big[\,\sU(F,t)\big].
$$
\end{lem}
\pf For all $w \notin  A$,  using the strong Markov property and the fact that   $F(X_s)=0$ for all $s<\sigma_{A}$, we get that 
\begin{align*}
	&\E^w\big[\, \sU(F, t) \big] 
			 = \E^w\big[ \, \sU(F, t) ; \sigma_{  A} \le t\big]
		= \E \bigg[ \E^{X_{\sigma_{A}}}\int_{0}^{t-\sigma_{A}} F(X_s) ds ;\, \sigma_{  A} \le t \bigg] 
		\le \sup_{x \in  A} \E^x \big[ \, \sU(F, t) \big].
\end{align*}
This proves the lemma. \qed

In the following, we give some consequences of our assumptions given in Section \ref{s:setting}. The following is a near infinity counterpart of \ref{A4}:

\bigskip

\textit{There exist constants $R_\infty \ge 1$, $\up \in (0,1)$,   $C_i>0$, $4\le i \le 7$  such that for every $x \in M$ and  $r>R_\infty \d(x)^\up$,}
\begin{equation}\label{B4}\tag{${\rm SP}^{R_\infty}(\varPhi,\up)$}
	C_4e^{-C_5 n}  \le \P^x\big(\tau_{B(x,r)} \ge n \varPhi(r)\big) \le  	C_6e^{-C_7 n}  \quad \text{for all} \; n \ge 1.
\end{equation}

\

By following the arguments in the proof of \cite[Proposition 4.3]{CKL} line-by-line, we obtain the next result. We skip the proof since the proof is almost identical. 

\begin{prop}
\label{p:E}
	(i)  Suppose that \VDo, \ref{A2} and \ref{A4+} hold for  an open set $U \subset M$. Then  \ref{A4} holds true with redefined $R_0>0$ and $C_0\in (0,1)$. Moreover, with the redefined $R_0$, \eqref{e:mean-exit} holds for every $x \in U$ and $0<r<R_0 \wedge (C_0 \updelta_U(x))$.
	
	\noindent	(ii) Suppose that \VDi, \ref{B2} and \ref{B4+} hold. Then  \ref{B4} holds  with a redefined $R_\infty\ge1$.  Moreover, with the redefined $R_\infty$, \eqref{e:mean-exit} holds for every $x \in M$ and $r>R_\infty \d(x)^\up$.
\end{prop}

We recall two propositions from \cite[Section 5]{CKL2} which will be used in this paper several times.

\begin{prop}%
	[\hspace{-0.1mm}{\cite[Proposition 5.1]{CKL2}}]\label{p:EP} 
	\noindent (i) Suppose that \ref{A1}, \ref{A2}, \ref{A3} and \ref{A4} hold for an open set $U\subset M$. Then there exist constants $q \in (0,1]$ and $c_1>0$ such that  for all $x \in U$, $0<r<3^{-1}( R_0 \wedge (C_0\updelta_U(x)))$ and $t>0$,
$$
		\P^x(\tau_{B(x,r)} \le t) \le c_1 \left(\frac{t}{\phi(x,r)}\right)^{q}.
$$
	
	\noindent (ii) Suppose that \VDi, \ref{B2}, \ref{B3} and \ref{B4+} hold. Then for every $\up_1 \in (\up, 1)$, there exist constants $c_1>0$ and  $R_1 \ge R_\infty$  such that 
	for all $x \in M$, $r>R_1 \d(x)^{\up_1}$ and $t \ge \vp(2r^{\up/\up_1})$,
		\begin{equation}\label{e:EP+n}
		\P^x(\tau_{B(x,r)} \le t) \le c_1 \frac{t}{\vp(r)}.
	\end{equation}
		Moreover, $X$ is conservative, that is, $\P^x(\zeta=\infty)=1$  for all $x \in M$.
\end{prop}

Let $(\theta_t)_{t\ge0}$ denote the shift operator which is defined by 
$(X_s \circ \theta_t) ( \omega) = X_{s+t}(\omega)$
 for all $t,s \ge 0$

An event $G$ is called \textit{shift-invariant} (with respect to $X$) if $G$ is a tail event, namely, $\cap_{t>0}^\infty \sigma(X_s:s>t)$-measurable,  and $\P^y(G)= \P^y(G \circ \theta_t)$ for all $y \in M$ and $t>0$.

\begin{prop}
[\hspace{-0.1mm}{\cite[Proposition 5.4]{CKL2}}]\label{p:law01}
	Suppose that \VDi, \ref{B2}, \ref{B3} and \ref{B4+} hold. Then for every shift-invariant event $G$, it holds either $\P^z(G)=0$ for all $z \in M$ or else $\P^z(G) = 1$ for all $z \in M$.
\end{prop}

We present limsup LILs for  $\tau_{B(x,r)}$ both at zero and at infinity. Note that the following result is closely related to Chung-type liminf LILs   \cite[Theorem 1.2 and Corollary 1.7]{CKL2}.

\begin{prop}\label{p:liminfLIL}
	(i)	Suppose that  \ref{A1}, \ref{A2}, \ref{A3} and \ref{A4} hold for an open set $U\subset M$.  Then, there exist constants $c_2 \ge c_1>0$  such that for every $x \in U$, there exists a constant $c_x \in [c_1,c_2]$ satisfying
	\begin{equation*}
		\limsup_{r\to0} \frac{\tau_{B(x,r)}}{\phi(x,r)\log|\log \phi(x,r)|}= c_x,~\quad\,
		\P^x\mbox{-a.s.}
	\end{equation*}
	
	\noindent(ii) 	Suppose that \VDi, \ref{B2}, \ref{B3} and \ref{B4+} hold.  Then, there exists a constant $c_3 \in (0, \infty)$ such that for every $x \in M$,
	\begin{equation}\label{e:limsup_tau}
		\limsup_{r\to\infty} \frac{\tau_{B(x,r)}}{\vp(r)\log\log \vp(r)} = c_3,~\quad\,
		\P^z\mbox{-a.s.}, \;\; \forall z \in M.
	\end{equation}
\end{prop}
\pf (i) By the Blumenthal's zero-one law, we get the result from  \cite[(5.11)]{CKL2}.

\noindent (ii)  Similarly, using Proposition \ref{p:law01}, we deduce the result from  \cite[(5.19)]{CKL2}. \qed

\begin{remark}
	{\rm In Proposition \ref{p:liminfLIL}(ii), we assumed \VDi \ and \ref{B4+}, which are stronger than \ref{B4} by Proposition \ref{p:E}(ii), to use  Proposition \ref{p:law01} and get a deterministic limit in \eqref{e:limsup_tau}.	
	}
\end{remark}

The following is one of  key estimates in the proof of  Theorems \ref{t:generalocc-SP} and \ref{t:generalocc2}.

\begin{lem}\label{l:mean-exponential}
	(i) Suppose that \ref{A1}, \ref{A2}, \ref{A3} and \ref{A4}  hold for an open set $U\subset M$.  Then there  exist constants $C_8,C_9,C_{10}>0$ such that for all $x  \in U$,  $0<r < 3^{-1}\big(R_0 \land (C_0 \delta_U(x))\big)$ and $y \in B(x,r)$,
	\begin{equation}\label{e:large_deviation_2}
		\E^y\bigg[\exp \Big(\,\frac{C_8\tau_{B(y,r)}}{ \phi(x,r)}\Big)\bigg] \le   e^{C_9}
	\end{equation}
	and
	\begin{equation}\label{e:large_deviation_1}
		\E^y\bigg[\exp \Big(-\frac{C_{10}\tau_{B(y,r/2)}}{ \phi(x,r)}\Big)\bigg] \le 2e^{-4}.
	\end{equation} 
	
	\noindent (ii) Suppose that \ref{B2}, \ref{B3} and \ref{B4}  hold. Define a function  $\vp$ by \eqref{e:molli1}. There exist constants $R_2 \ge R_\infty$ and $C_8,C_9,C_{10}>0$ such that  \eqref{e:large_deviation_2} and \eqref{e:large_deviation_1} hold for all $x \in M$, $r \ge R_2 \d(x)^{\up^{2/3}}$ and  $y \in  B(x, r^{\up^{-1/3}})$, with  $\vp(r)$ in denominators instead of $\phi(x,r)$.
\end{lem}
\pf (i)  Choose any $x \in U$ and let $r_1:= R_0 \wedge (C_0\updelta_U(x))$. By the triangle inequality, we see that for all $0<r<r_1/3$ and $y \in B(x,r)$,
\begin{equation*}
	 R_0 \wedge (C_0\updelta_U(y)) \ge  R_0 \wedge (C_0\updelta_U(x)-C_0 r) \ge r_1/2.
\end{equation*}
Hence, we get from \ref{A1} and \ref{A4} that for all $0<r<r_1/3$ and $y \in B(x,r)$,  
\begin{align}\label{e:exponential1}
	\E^y\left[\exp \Big(\frac{C_7}{2C_1}\frac{\tau_{B(y,r)}}{ \phi(x,r)}\Big)\right] &\le \E^y\left[\exp \Big(\frac{C_7}{2}\frac{\tau_{B(y,r)}}{ \phi(y,r)}\Big)\right] \le \sum_{n=1}^\infty  e^{C_7n/2} \, \P^y\big( (n-1) \phi(y,r) \le \tau_{B(y,r)} < n \phi(y,r)\big)\nn\\
	& \le  \sum_{n=1}^\infty  e^{C_7n/2} \, \P^y\big(  \tau_{B(y,r)} \ge (n-1) \phi(y,r)\big) \le  C_6 e^{C_7}\sum_{n=1}^\infty e^{-C_7n/2 }<\infty.
\end{align}
Moreover, by \ref{A1}, \ref{A2} and Proposition \ref{p:EP}(i),  there exists a large constant  $c_1>0$  independent of $x$ such that for all $0<r<3^{-1}(R_0 \land (C_0 \delta_U(x)))$ and $y \in B(x,r)$,
\begin{equation}\label{e:SP_EP}
	\P^y\big(\tau_{B(y,r/2)}\le 4c_1^{-1}\phi(x,r)\big)  \le 	\P^y\big(\tau_{B(y,r/2)}\le 4c_1^{-1}C_1C_2\phi(y,r/2)\big) \le  e^{-4}.
\end{equation} 
Then using Markov inequality, we obtain that  for all $0<r<3^{-1}(R_0 \land (C_0 \delta_U(x)))$ and $y \in B(x,r)$,
\begin{align}\label{e:exponential2}
	&\E^y\left[\exp \Big(-\frac{c_1\tau_{B(y,r/2)}}{ \phi(x, r)}\Big)\right]\nn\\
	&\le \P^y\big( \tau_{B(y,r/2)} \le  4c_1^{-1} \phi(x,r) \big) + \E^y\left[\exp \Big(-\frac{c_1\tau_{B(y,r/2)}}{ \phi(x,r)}\Big)  :   \frac{c_1\tau_{B(y,r/2)}}{\phi(x,r)}> 4 \right]\le 2 e^{-4}. 
\end{align}

\noindent (ii)  Choose any $x \in M$ and set $\up_1:= \up^{2/3} \in (\up,1)$. Let $R_1 \ge R_\infty$ be the constant from Proposition \ref{p:EP}(ii). 
Note that for all  $r \ge  (2R_1)^{1/(1-\up^{1/3})}\d(x)^{\up_1}$, we have
$ r > 2R_1 \d(x)^{\up_1}$ and $r^{1-\up^{1/3}} \ge 2R_1$ since  $R_1 \ge 1$ and $\d(x) \ge 1$. Therefore,  for all $r \ge  (2R_1)^{1/(1-\up^{1/3})}\d(x)^{\up_1}$ and $y \in B(x, r^{\up^{-1/3}})$,
\begin{equation*}
	R_1\d(y)^{\up_1} \le R_1(\d(x)+ d(x,y))^{\up_1}  \le R_1 \d(x)^{\up_1} + R_1r^{\up_1 \up^{-1/3}} = R_1 \d(x)^{\up_1} + R_1r^{\up^{1/3}}< \frac{r}{2} + \frac{r}{2}=r.
\end{equation*}
Hence, using    \eqref{e:molli2}, \ref{B2}, \ref{B4}  and Proposition \ref{p:EP}(ii), by similar arguments as that for (i),  we can deduce that \eqref{e:exponential1} and \eqref{e:exponential2} hold for all  $r \ge  (2R_1)^{1/(1-\up^{1/3})}\d(x)^{\up_1}$ and $y \in B(x, r^{\up_1/\up})$ with $\vp(r)$ in denominators instead of $\phi(x,r)$. We omit the details. \qed

\section{Proofs of Theorems \ref{t:generalocc-SP}, \ref{t:generalocc2} and \ref{t:generalocc3}}\label{s:occ}

In this section, we give proofs for the first type of LILs: Theorems \ref{t:generalocc-SP}, \ref{t:generalocc2} and  \ref{t:generalocc3}. 
The following lemma shows the monotone properties of the functions $f_0(x, \cdot)$ and $f_\infty(\cdot)$ in Theorems \ref{t:generalocc-SP} and \ref{t:generalocc2}.

\begin{lem}\label{l:occ-lem-1}
Let $f$ and $g$ be increasing positive continuous   functions defined on a subinterval of $(0,\infty)$ such that $f \le g$.

\smallskip

\noindent(i) If $\lim_{r \to 0}g(r)=0$, then for every $x \in M$,
\begin{equation*}
	\limsup_{r \to 0} \frac{\sU(B(x,r), f(r))}{f(r)} \ge 	\limsup_{r \to 0} \frac{\sU(B(x,r), g(r))}{g(r)}.
\end{equation*}
In particular, for each fixed $x \in U$, the function $f_0(x, \cdot)$ defined as \eqref{e:occ_SP} is non-increasing.

\noindent(ii) If $\lim_{r \to \infty}f(r)=\infty$, then for every $x \in M$,
\begin{equation*}
\limsup_{r \to \infty} \frac{\sU(B(x,r), f(r))}{f(r)} \ge 	\limsup_{r \to \infty} \frac{\sU(B(x,r), g(r))}{g(r)}.
\end{equation*}
In particular, if the function $f_\infty$ in \eqref{e:generalocc2_1} is well-defined, then it is non-increasing.
\end{lem}
\pf (i) Note that  $(f^{-1} \circ g)(r) \ge (f^{-1} \circ f)(r) =r$ for all $r>0$.  Hence by the monotonicity of occupation times, we get
\begin{align*}
	\limsup_{r \to 0} \frac{\sU(B(x,r), f(r))}{f(r)} = \limsup_{r \to 0} \frac{\sU(B(x, ( f^{-1} \circ g)(r)), g(r))}{g(r)} \ge 	\limsup_{r \to 0} \frac{\sU(B(x,r), g(r))}{g(r)}.
\end{align*}

\noindent (ii) It can be proved by the same way. \qed

Next, we prove the property (P2) in Theorem \ref{t:generalocc-SP} and the first equality in \eqref{e:finfty}.

\begin{lem}\label{l:occ_lower}
	(i) Under the setting of Theorem \ref{t:generalocc-SP},  there exists a constant $\kappa_0>0$  such that for all $x \in U$ and $\kappa \le \kappa_0$,
	\begin{equation}\label{e:occ_lower_zero}
		\limsup_{r \to 0} \frac{\sU\big(B(x, r),\, \kappa \phi(x,r) \log |\log \phi(x,r)|\big)}{\kappa\phi(x,r) \log |\log \phi(x,r)|} =1,~\quad\,
		\P^x\mbox{-a.s.}
	\end{equation}
	
	\noindent(ii) Under the setting of Theorem \ref{t:generalocc2}, there exists a constant $\kappa_\infty>0$ such that for all $x \in M$ and $\kappa \le \kappa_\infty$,
	\begin{equation}\label{e:occ_lower_inf}
		\limsup_{r \to \infty} \frac{\sU\big(B(x, r),\, \kappa \vp(r) \log \log \vp(r)\big)}{\kappa\vp(r) \log \log \vp(r)} =1,~\quad\,
		\P^z\mbox{-a.s.}, \;\; \forall z \in M.
	\end{equation}
\end{lem}
\pf (i)  By Proposition  \ref{p:liminfLIL}(i), there exists a constant $\kappa_0>0$ independent of $x \in U$ such that for $\P^x$-a.s $\omega$, there is a decreasing sequence $(r_n)_{n \ge 1}=(r_n(\omega))_{n \ge 1}$ 
which converges to zero and
\begin{equation*}
	\tau_{B(x, r_n)} \ge \kappa_0 \phi(x,r_n) \log |\log \phi(x,r_n)| \quad \text{for all} \;\; n \ge 1.
\end{equation*}
Then for all $\kappa\le \kappa_0$ and  $n \ge 1$, it holds that 
\begin{align*}
		\sU\big(B(x,r_n), \kappa \phi(x,r_n) \log |\log  \phi(x,r_n)|\big)(\omega) &\ge \tau_{B(x,r_n)}(\omega) \wedge \big( \kappa \phi(x,r_n) \log |\log  \phi(x,r_n)| \big)\\
	&=\kappa \phi(x,r_n) \log |\log  \phi(x,r_n)| .
\end{align*}
This proves \eqref{e:occ_lower_zero}.

\noindent (ii) Analogously, one can deduce \eqref{e:occ_lower_inf} from Proposition \ref{p:liminfLIL}(ii).  \qed

In the following, we show that the function $f_0$ defined as \eqref{e:occ_SP} satisfies that $\lim_{\kappa \to \infty} f_0(x, \kappa)=0$ uniformly on $U$. In order to prove this, we give a number of definitions first.

  Let $x \in U$.   For $\kappa>0$ and $\delta \in (0,1/4)$, let
$(u_n^{\kappa, \delta})_{n \ge N}=(u_n^{\kappa, \delta}(x))_{n \ge N}$ be a decreasing sequence  such that 
\begin{equation}\label{e:def_u_n}
	\kappa \phi(x,u_n^{\kappa, \delta})\log|\log\phi(x,u_n^{\kappa, \delta} )|=(1-\delta)^n  \quad \text{for all} \;\; n\ge N.
\end{equation}
By taking logarithm in \eqref{e:def_u_n} twice and assuming that $N$ is large enough, we get that
\begin{equation}\label{e:loglog_u_n}
	2^{-1}\log n \le \log |\log \phi(x,u^{\kappa, \delta}_n)| \le 2\log n  \quad \text{for all} \;\; n\ge N.
\end{equation}
For $n,m \ge 1$, $\kappa>0$ and $\delta \in (0,1/4)$, define an event $ O_n(\kappa,\delta,m;x)$ as
\begin{equation}\label{e:def_O_n}
	O_n(\kappa,\delta,m;x)= 
	\left\{\, \omega \in \Omega\,;\, \sU \big(B(x,u_n^{\kappa, \delta}),(1-\delta)^n\big)(\omega) \ge (1-\delta)^{n+1}(1-m\delta)\right\}.
\end{equation}
Lastly, we  define a constant $L$ which is  independent of $x \in U$ as 
\begin{equation}\label{e:def_L}
	L= 4^{1+\beta_2}C_U  (2C_9+1)C_{10}/C_8,
\end{equation}
where $\beta_2,C_U,C_8, C_9, C_{10}$ are the positive constants from  \eqref{e:scaling-zero}, \eqref{e:large_deviation_2} and  \eqref{e:large_deviation_1}.

Now we claim that there exists a constant $\delta_0 \in (0, 1/4)$ such that for every $x \in U$, $\delta \in (0,\delta_0]$ and $m \in \N$ satisfying $ (m + L) \delta \le 1 $, there exists a constant $\kappa_{\delta,m} > 0$ independent of $x$ such that 
\begin{equation}\label{e:Oni}
	\sum_{n=1}^\infty \P^x\big(O_n(\kappa_{\delta,m},\delta,m;x)\big) < \infty.
\end{equation}

Before giving the proof of  \eqref{e:Oni}, we note that \eqref{e:Oni} yields the desired result.

\begin{lem}\label{l:generalocc-SP-4}
	If \eqref{e:Oni} is true, then $\lim_{\kappa \to \infty}f_0(x, \kappa)=0$ uniformly on $U$.
\end{lem}
\pf  Choose any $x \in U$ and $\delta \in (0, \delta_0/(L+1))$.  Let $m:=\lfloor \delta^{-1} -L \rfloor$ and $\lambda:=\kappa_{\delta,m}$. By the Borel-Cantelli lemma, it holds $\P^x$-a.s. that 
\begin{align}\label{e:generalocc_last}
	&\limsup_{r \to 0} \frac{\sU \big(B(x, r), \lambda\phi(x,r)\log |\log \phi(x,r)|\big)}{\lambda\phi(x,r)\log |\log \phi(x,r)|} \nn\\
	&\le \limsup_{n \to \infty} \; \sup \left \{ \frac{\sU \big(B(x, r),\lambda \phi(x,r)\log |\log \phi(x,r)|\big)}{\lambda\phi(x,r)\log |\log \phi(x,r)|} : u^{\lambda, \delta }_{n+1}\le r < u^{\lambda,\delta }_n  \right\}\nn\\
	& \le \limsup_{n \to \infty} \frac{\sU  \big(B(x,  u_n^{\lambda, \delta}), \lambda\phi(x,u_n^{\lambda, \delta})\log|\log \phi(x,u_n^{\lambda, \delta})|\big)}{\lambda\phi(x,u_{n+1}^{\lambda, \delta})\log|\log \phi(x,u_{n+1}^{\lambda, \delta})|}=\limsup_{n \to \infty} \frac{\sU \big(B(x,  u_n^{\lambda, \delta}),  (1-\delta)^n\big)}{(1-\delta)^{n+1}}\nn\\
	&\le \limsup_{n \to \infty} \frac{(1-\delta)^{n+1}(1-m \delta)}{(1-\delta)^{n+1}} = 1- m \delta \le (L+1)\delta.
\end{align}
Hence, $f_0(x, \kappa) \le (L+1)\delta$ for all $\kappa \ge \lambda$ by Lemma \ref{l:occ-lem-1}. Recall that the constant $\lambda=\kappa_{\delta,m}$ is independent of $x$.  Since $\delta$ can be arbitrarily small, the proof is finished.\qed 

In the following two lemmas, we  prove \eqref{e:Oni} by induction. Recall that $C_i$, $i=8,9,10$ are the  constants from \eqref{e:large_deviation_2} and \eqref{e:large_deviation_1}.

\begin{lem}\label{l:generalocc-SP-2}
	Under the setting of Theorem \ref{t:generalocc-SP}, there exist constants $\delta_0 \in (0,1/4)$ and $\lambda_1>0$ such that for any $x \in U$ and  $\delta \in (0,\delta_0]$,  \eqref{e:Oni} holds when $m=1$ and $\kappa_{\delta,1}=\lambda_1$. 
\end{lem}
\pf Let $\delta_0 \in (0,1/4)$ and $\lambda_1>0$ be constants chosen later. Fix $x \in U$ and $\delta \in (0, \delta_0]$. We simply write  $u_n=u_n^{\lambda_1, \delta}(x)$ and $O_n=O_n(\lambda_1, \delta,1;x)$.  For each $n \ge N$, we  define sequences of  random times $(\sigma^n_j)_{j \ge 1}$ and $(\tau^n_j)_{j \ge 1}$ as
\begin{equation}\label{e:out-and-in}
	\sigma_1^n=0, \quad 	\tau^n_{j}= \inf\{s >\sigma^n_{j}: X_s \notin  B(x, 2u_n)  \} \quad \text{and} \quad \sigma^n_{j+1}= \inf\{s >\tau^n_{j}: X_{s} \in B(x,  u_n)  \}.
\end{equation}

Recall that $(\theta_t)_{t\ge 0}$ denotes the shift operator. By Markov inequality, we have that 
\begin{align}\label{e:On-tau-le-log}
	&\P^x\left(O_n  :  \tau^n_{\lfloor \log n \rfloor} \ge (1-\delta)^n\right) 	\le \P^x 
	\Big(O_n : \sum_{j=1}^{\lfloor \log n \rfloor}\tau_{B(x, 2 u_n)} \circ \theta_{\sigma^n_j}  \ge \sU\big(B(x, u_n), (1-\delta)^n\big) \Big)\nn\\
	&\le \P^x 
	\bigg(\sum_{j=1}^{\lfloor \log n \rfloor}\tau_{B(x, 2 u_n)} \circ \theta_{\sigma^n_j}  \ge  (1-\delta)^{n+2}  \bigg)\nn\\
	&\le  \exp \Big(- \frac{C_8(1-\delta)^{n+2}}{\phi(x,4  u_n)} \Big) \, \E^x\bigg[\exp \Big(\frac{C_8}{\phi(x,4u_n)} \sum_{j=1}^{\lfloor \log n \rfloor}\tau_{B(x, 2 u_n)} \circ \theta_{\sigma^n_j}\Big) \bigg].
\end{align}
Using the strong Markov property and the fact that $B(x, 2u_n) \subset B(y, 4u_n)$ for all $y \in B(x, 2 u_n)$ in the first inequality below and  \eqref{e:large_deviation_2} in the last inequality, we get that for all $n$ large enough, 
\begin{align}\label{e:generalocc-eq1}
	&\E^x\bigg[\exp \Big(\frac{C_8}{\phi(x,4u_n)} \sum_{j=1}^{\lfloor \log n \rfloor}\tau_{B(x, 2 u_n)} \circ \theta_{\sigma^n_j}\Big) \bigg] \nn\\
	&= \E\E^x\bigg[\exp \Big(\frac{C_8}{\phi(x,4u_n)} \sum_{j=1}^{\lfloor \log n \rfloor}\tau_{B(x, 2 u_n)} \circ \theta_{\sigma^n_j}\Big)\,\Big|\, \FF_{\sigma^n_{\lfloor \log n \rfloor}}  \bigg]\nn\\
	&\le  \E^x\bigg[\exp \Big(\frac{C_8}{\phi(x,4u_n)} \sum_{j=1}^{\lfloor \log n \rfloor-1}\tau_{B(x, 2 u_n)} \circ \theta_{\sigma^n_j}\Big) \bigg] \sup_{y \in  B(x, 2 u_n)}  \E^y\bigg[\exp \Big(\frac{C_8\tau_{B(y, 4u_n)}}{\phi(x,4u_n)} \Big) \bigg]\nn \\
	& \le \cdots \le  \bigg(\sup_{y \in  B(x, 2 u_n)} \E^y\bigg[\exp \Big(\frac{C_8  \tau_{B(y, 4u_n)}}{\phi(x,4u_n)} \Big) \bigg] \bigg)^{\lfloor \log n \rfloor }\le e^{C_9 \log n}.
\end{align}
By \eqref{e:On-tau-le-log}, \eqref{e:generalocc-eq1} and the definition \eqref{e:def_u_n}, using \eqref{e:loglog_u_n},   \eqref{e:scaling-zero} and   the fact that $\delta<1/4$, we get that for all $n$ large enough, 
\begin{align}\label{e:On_small_Tn}
\P^x\left(O_n  :  \tau^n_{\lfloor \log n \rfloor} \ge (1-\delta)^n\right) &\le  \exp \Big( C_9 \log n - \lambda_1 \frac{ C_8(1-\delta)^{2} \phi(x,u_n) \log |\log \phi(x,u_n)| }{\phi(x,4  u_n)} \Big) \nn\\
	&\le  \exp\Big(\big(C_9- c_1\lambda_1\big)\log n\Big).
\end{align}
where $c_1:=4^{-1-\beta_2}C_8 / C_U$ is independent of  $(x,\kappa_1, \delta)$.

Set $r_0:=3^{-1}(R_0 \wedge (C_0\updelta_U(x)))$. Using  Markov inequality and the fact that  $\sU(D_1,t) \le t-\sU(D_2,t)$ for all $t \ge 0$ when $D_1, D_2 \in \sB(M)$ are disjoint, we see that for all $n\ge N$,
\begin{align}\label{e:generalocc-eq2}
	&\P^x\left(O_n  :  \tau^n_{\lfloor \log n \rfloor} < (1-\delta)^n < \tau_{B(x,r_0)} \right) \nn\\
	&\le \P^x 
	\bigg(O_n : \sum_{j=1}^{\lfloor \log n \rfloor-1}\tau_{B(X_0, d(x,X_0)-u_n)} \circ \theta_{\tau^n_j}   \le (1-\delta)^n - \sU\big(B(x, u_n),   (1-\delta)^n\big),  \nn\\
	& \hspace{25mm} \tau^n_{\lfloor \log n \rfloor} < (1-\delta)^n< \tau_{B(x,r_0)} \bigg) \nn\\
	&\le \P^x 
	\bigg(\sum_{j=1}^{\lfloor \log n \rfloor-1}\tau_{B(X_0, d(x,X_0)-u_n)} \circ \theta_{\tau^n_j}  \le 2\delta(1-\delta)^n, \;  X({\tau_j^n}) \in B(x,r_0),\, 1\le j \le \lfloor \log n \rfloor-1 \bigg) \nn\\
	&\le  \exp \Big( \frac{2C_{10}\delta (1-\delta)^n}{\phi(x,2u_n)} \Big)  \, \E^x\bigg[\exp \Big(-\frac{C_{10}}{\phi(x,2u_n)} \sum_{j=1}^{\lfloor \log n \rfloor-1}\tau_{B(X_0, d(x,X_0)-u_n)} \circ \theta_{\tau^n_j}\Big)\nn\\
	& \hspace{68mm} :    X({\tau_j^n}) \in B(x,r_0),\, 1\le j \le \lfloor \log n \rfloor-1\bigg].
\end{align}
For every $y \in B(x,r_0) \setminus B(x,2u_n)$,  we get from \eqref{e:large_deviation_1} that
\begin{equation}\label{e:generalocc-eq3-0}
	\E^x\bigg[\exp \Big(-\frac{C_{10}}{\phi(x,2u_n)} \tau_{B(y, d(x,y)-u_n)}\Big)\bigg] \le 	\E^x\bigg[\exp \Big(-\frac{C_{10}}{\phi(x,d(x,y))} \tau_{B(y, d(x,y)/2)}\Big)\bigg]\le e^{-3}.
\end{equation}
Since $X(\tau_j^n) \notin B(x,2u_n)$, by  the strong Markov property, it follows that for all $n$ large enough, 
\begin{align}\label{e:generalocc-eq3}
	&\E^x\bigg[\exp \Big(-\frac{C_{10}}{\phi(x,2u_n)} \sum_{j=1}^{\lfloor \log n \rfloor-1}\tau_{B(X_0, d(x,X_0)-u_n)} \circ \theta_{\tau^n_j}\Big) :    X({\tau_j^n}) \in B(x,r_0),\, 1\le j \le \lfloor \log n \rfloor-1\bigg] \nn\\
	& = \E\E^x\bigg[\exp \Big(-\frac{C_{10}}{\phi(x,2u_n)} \sum_{j=1}^{\lfloor \log n \rfloor-1}\tau_{B(X_0, d(x,X_0)-u_n)} \circ \theta_{\tau^n_j}\Big) \nn\\
	& \hspace{1.5cm} :  X({\tau_j^n}) \in B(x,r_0),\, 1\le j \le \lfloor \log n \rfloor-1\,\Big|\, \FF_{\tau^n_{\lfloor \log n \rfloor-1}}  \bigg] \nn\\
	&\le  \E^x\bigg[\exp \Big(-\frac{C_{10}}{\phi(x,2u_n)} \sum_{j=1}^{\lfloor  \log n \rfloor-2}\tau_{B(X_0, d(x,X_0)-u_n)} \circ \theta_{\tau^n_j}\Big)  :  X({\tau_j^n}) \in B(x,r_0),\, 1\le j \le \lfloor \log n \rfloor-2  \bigg] \nn\\ 
	& \quad\quad \times \sup_{y \in  B(x, r_0) \setminus B(x,2u_n)}  \E^y\bigg[\exp \Big(-\frac{C_{10}\tau_{B(y, d(x,y) -u_n)}}{\phi(x,2u_n)} \Big) \bigg] \nn\\
	& \le \cdots \le  \bigg(  \sup_{y \in  B(x, r_0) \setminus B(x,2u_n)}  \E^y\bigg[\exp \Big(-\frac{C_{10}\tau_{B(y, d(x,y) -u_n)}}{\phi(x,2u_n)} \Big) \bigg] \bigg)^{\lfloor \log n \rfloor -1} \hskip -0.05in \le e^{6- 3\log n}.
\end{align}
Therefore, using the definition \eqref{e:def_u_n}, the monotone property of $\phi(x, \cdot)$ and \eqref{e:loglog_u_n}, we obtain that for all $n$ large enough, 
\begin{align}\label{e:On_large_Tn}
\P^x\left(O_n  :  \tau^n_{\lfloor \log n \rfloor} < (1-\delta)^n < \tau_{B(x,r_0)} \right) 
	&\le e^6 \exp \Big(  -3 \log n + \delta \lambda_1 \frac{2C_{10}  \phi(x,u_n) \log |\log \phi(x,u_n)|}{\phi(x,2 u_n)} \Big)\nn\\
	&\le e^6 \exp\Big(\big(-3+4C_{10}\delta \lambda_1\big)\log n\Big).
\end{align}

Take  $\lambda_1= (2+C_9)/c_1$ and  $\delta_0=1/ (4C_{10} \lambda_1)$. Then by \eqref{e:On_small_Tn},  \eqref{e:On_large_Tn} and Proposition  \ref{p:EP}(i), it holds that for all $n$ large enough,
\begin{align*}
	&\P^x\big(O_n (\lambda_1, \delta, 1; x)\big) \nn\\[5pt]
	&\le \P^x\big(O_n  :  \tau^n_{\lfloor \log n \rfloor} \ge (1-\delta)^n \big) + \P^x\big(O_n  :  \tau^n_{\lfloor \log n \rfloor} < (1-\delta)^n< \tau_{B(x,r_0)} \big) + \P^x\big( \tau_{B(x,r_0)} \le (1-\delta)^n\big)\\[4pt]
	& \le \exp\Big((C_9-c_1\lambda_1)\log n\Big) + e^6\exp\Big((-3+4C_{10}\delta_0 \lambda_1)\log n\Big) +c_2 \phi(x,r_0)^{-q} (1-\delta)^{nq}\\[4pt]
	& \le n^{-2} + e^6 n^{-2} + c_2 \phi(x,r_0)^{-q} (1-\delta)^{nq}.
\end{align*} 
The proof is complete. \qed

\begin{lem}\label{l:generalocc-SP-3}
	Under the setting of Theorem \ref{t:generalocc-SP},  \eqref{e:Oni} is true.
\end{lem}
\pf Choose any $x \in U$ and $\delta \in (0, \delta_0]$. When $m =1$, \eqref{e:Oni} is satisfied by Lemma \ref{l:generalocc-SP-2}. For the next step, let $m \in \N$ be such that  $(m+1+L)\delta  \le 1$ and suppose that \eqref{e:Oni} is satisfied with $m$. 

 Let $\kappa_{\delta,m+1} \ge 4^{1+\beta_2}C_U\kappa_{\delta,m}$ be a constant  chosen later. For simplicity, we write as 
$$
u_n = u_n^{\kappa_{\delta, m},\,\delta}(x), \;\; \wt u_n = u_n^{\kappa_{\delta, m+1},\,\delta}(x), \;\; O_n=O_n(\kappa_{\delta, m},\delta,m;x)  \;\; \text{ and } \;\;   \wt O_n=O_n(\kappa_{\delta, m+1},\delta,m+1;x),
$$
where $u_n^{\kappa, \, \delta}(x)$ and $O_n(\kappa, \delta,m;x)$ are defined by \eqref{e:def_u_n} and \eqref{e:def_O_n} respectively. 
Using \eqref{e:scaling-zero} and  \eqref{e:loglog_u_n}, we get that  for all $n$ large enough,
\begin{align}\label{e:wtu_n_u_n}
		4\kappa_{\delta, m}\phi(x,4\wt u_n)& \le 	4^{1+\beta_2 }C_U\kappa_{\delta,m} \phi(x,\wt u_n) \le \kappa_{\delta, m+1} \phi(x,\wt u_n)= \frac{(1-\delta)^n}{\log|\log \phi(x,\wt u_n)|} \nn\\
	&\le \frac{2(1-\delta)^n}{\log n}  \le  \frac{4(1-\delta)^n}{\log|\log \phi(x,u_n)|} = 4\kappa_{\delta, m} \phi(x,u_n).
\end{align}
Thus, $4 \wt u_n \le u_n$ for all $n$ large enough. 

   Our task is to show $\sum_{n=1}^\infty \P^x(\wt O_n)<\infty$. Since $\sum_{n=1}^\infty \P^x(O_n) < \infty$ by the induction hypothesis, it suffices to show  $\sum_{n=1}^\infty \P^x(\wt O_n \setminus O_n) < \infty$. Define $(\wt \sigma^n_j)_{j \ge 1}$ and $(\wt \tau^n_j)_{j \ge 1}$ as  (cf. \eqref{e:out-and-in})
\begin{align*}
	\wt \sigma_1^n:=0, \;\; \wt 	\tau^n_{j}:= \inf\{s >\wt \sigma^n_{j}: X_s \notin  B(x, 2\wt u_n)  \} \quad \text{and} \quad \wt  \sigma^n_{j+1}:= \inf\{s >\wt \tau^n_{j}: X_{s} \in B(x,  \wt u_n)  \}.
\end{align*}
Following the calculations in \eqref{e:On-tau-le-log}, we see that for any constant $K>0$,
\begin{align*}
	&\P^x\Big(\wt O_n \setminus O_n  :  \wt \tau^n_{\,2\lfloor K\log n \rfloor} \ge (1-\delta)^n\Big)\\
	& 	\le \P^x 
	\Big(\wt O_n \setminus O_n : \sum_{j=1}^{2\lfloor K\log n \rfloor}  \tau_{B(x,2\wt u_n)} \circ \theta_{\wt \sigma^n_j}  \ge \sU\big(B(x, \wt u_n), (1-\delta)^n\big) \Big)\nn\\
	&\le \P^x 
	\Big(\sum_{j=1}^{2\lfloor K\log n \rfloor}  \tau_{B(x, 2\wt u_n)} \circ \theta_{\wt \sigma^n_j}  \ge  (1-\delta)^{n+1} (1- \delta- m\delta)  \Big)\nn\\
	&\le  \exp \Big(- \frac{C_8(1-\delta)^{n+1}(1-\delta-  m\delta)}{\phi(x,4 \wt u_n)} \Big)\E^x\bigg[\exp \Big(\frac{C_8}{\phi(x,4\wt u_n)} \sum_{j=1}^{2\lfloor K\log n \rfloor}\tau_{B(x, 2 \wt u_n)} \circ \theta_{\wt \sigma^n_j}\Big) \bigg].
\end{align*}
Also, as in the calculations showing \eqref{e:generalocc-eq1}, we get from the strong Markov property and \eqref{e:large_deviation_2} that for all $n$ large enough, 
\begin{align*}
	&\E^x\bigg[\exp \Big(\frac{C_8}{\phi(x,4\wt u_n)} \sum_{j=1}^{2\lfloor K\log n \rfloor}  \tau_{B(x, 2 \wt u_n)} \circ \theta_{\wt \sigma^n_j}\Big) \bigg] \\
	&\le  \E^x\bigg[\exp \Big(\frac{C_8}{\phi(x,4 \wt u_n)} \sum_{j=1}^{2\lfloor K\log n \rfloor-1} \tau_{B(x, 2\wt u_n)} \circ \theta_{\wt \sigma^n_j}\Big) \Big] \sup_{y \in  B(x, 2 \wt u_n)}  \E^y\bigg[\exp \Big(\frac{C_8\tau_{B(y, 4\wt u_n)}}{\phi(x,4 \wt u_n)} \Big) \bigg] \\
	& \le \cdots \le  \bigg(\sup_{y \in  B(x, 2\wt u_n)} \E^y\bigg[\exp \Big(\frac{C_8 \tau_{B(y, 4\wt u_n)}}{\phi(x,4\wt u_n)} \Big) \bigg] \bigg)^{2\lfloor K\log n \rfloor }\le e^{2C_9K\log n}. 
\end{align*}
Thus, since $\delta\le \delta_0< 1/4$ and $1-\delta- m\delta \ge L\delta$, 
by \eqref{e:loglog_u_n} and \eqref{e:scaling-zero}, it holds that  for all $n$ large enough,
\begin{align}\label{4.1.1}
&	\P^x\Big(\wt O_n\setminus O_n  :  \wt \tau^n_{\,2\lfloor K\log n \rfloor}\ge (1-\delta)^n\Big)\nn\\[4pt]
&\le \exp \bigg( 2C_9 K \log n   - \kappa_{\delta,m+1}  \frac{C_8(1-\delta)(1-\delta- m\delta)  \phi(x,\wt u_n) \log |\log \phi(x,\wt u_n)|}{\phi(x,4\wt u_n)}   \bigg)  \nn \\[4pt]
	&  \le  \exp \Big( \big(2C_9 K - c_1 L \delta  \kappa_{\delta,m+1}\big)\log n\Big),
\end{align}
where  $c_1:=4^{-1-\beta_2}C_8/C_U$.

Next, we  set  $r_0:=3^{-1}(R_0 \wedge (C_0 \updelta_U(x)))$ as in the proof of Lemma \ref{l:generalocc-SP-2}, and define two random sequences of natural numbers  $(a_\ell)_{\ell \ge 0}$ and $(b_\ell)_{\ell \ge 0}$ as
\begin{equation}\label{e:def-a_n}
a_0=0, \quad 	a_{\ell +1}  := \inf\{  j > a_\ell : X({\wt \tau_j^n}) \in B(x,u_n - \wt u_n) \setminus B(x,2\wt u_n) \}
\end{equation}
and
\begin{equation}\label{e:def-b_n}
b_0=0, \quad 	b_{\ell +1}:= \inf\{ j > b_\ell : X({\wt \tau_j^n}) \notin B(x,u_n - \wt u_n)  \}.
\end{equation}
Let $\wt\tau^{n,1}_\ell:= \wt\tau^n_{a_\ell}$ and $\wt\tau^{n,2}_\ell := \wt \tau^n_{b_\ell}$. Since $X(\wt \tau_j^n) \notin B(x, 2 \wt u_n)$ by the definition of $\wt \tau_j^n$, we have
\begin{align}\label{4.1.2}
&	\P^x\Big(\wt O_n\setminus O_n  :  \wt \tau^n_{\,2\lfloor K\log n \rfloor} < (1-\delta)^n < \tau_{B(x,r_0)} \Big) \nn\\
	&\le 	\P^x\Big(\wt O_n\setminus O_n  :  \wt \tau^{n,1}_{\lfloor K\log n \rfloor} < (1-\delta)^n< \tau_{B(x,r_0)}\Big)\nn\\
	&\qquad \qquad  +	\P^x\Big(\wt O_n\setminus O_n  :  \wt \tau^{n,2}_{\lfloor K\log n \rfloor} < (1-\delta)^n< \tau_{B(x,r_0)}\Big). 
\end{align}
 
 Observe that for all $y \in B(x,u_n - \wt u_n) \setminus B(x,2\wt u_n)$ and  $\omega \in \wt O_n \setminus O_n$, 
\begin{align}\label{e:O_n_minus_O_n}
\sU\big(B(y,\wt u_n), (1-\delta)^n\big)(\omega) &\le 	\sU\big(B(x,u_n) \setminus B(x,\wt u_n) , (1-\delta)^n\big) (\omega) \nn\\
&= \sU\big(B(x,u_n), (1-\delta)^n) \big)(\omega) -  \sU\big(B(x,\wt u_n), (1-\delta)^n\big)(\omega) \nn\\
	&\le (1-\delta)^{n+1} \big( (1-m\delta) - (1-(m+1)\delta) \big)= \delta(1-\delta)^{n+1}. 
\end{align}
Thus, by following the calculations in  \eqref{e:generalocc-eq2}, \eqref{e:generalocc-eq3-0}, \eqref{e:generalocc-eq3} and \eqref{e:On_large_Tn} in turn, we get that for all $n$ large enough,
\begin{align}\label{4.1.2-1}
	&\P^x\Big(\wt O_n \setminus O_n:  \wt \tau^{n,1}_{\lfloor K\log n \rfloor} < (1-\delta)^n < \tau_{B(x,r_0)}\Big) \nn\\
	&\le \P^x 
	\bigg(\sum_{j=1}^{\lfloor K\log n \rfloor-1}\tau_{B(X_0, d(x,X_0) - \wt u_n)} \circ \theta_{\wt \tau^{n,1}_j}  \le \delta(1-\delta)^{n+1}, \;  X({\wt \tau_j^{n,1}}) \in B(x,r_0),\, 1\le j \le \lfloor K\log n \rfloor-1 \bigg) \nn\\
	&\le \exp \Big( \frac{C_{10}\delta(1-\delta)^{n+1}}{\phi(x,2 \wt u_n)} \Big) \E^x\bigg[\exp \Big(-\frac{C_{10}}{\phi(x,2\wt u_n)} \sum_{j=1}^{\lfloor K\log n \rfloor-1}\tau_{B(X_0,d(x,X_0) - \wt u_n)} \circ \theta_{\wt\tau^{n,1}_j}\Big)  \nn\\
	& \hspace{65mm} :    X({\wt \tau_j^{n,1}}) \in B(x,r_0),\, 1\le j \le \lfloor K\log n \rfloor-1\bigg] \nn\\
	& \le   \bigg(\sup_{y \in  B(x, r_0) \setminus B(x, 2 \wt u_n)}  \E^y\bigg[\exp \Big(-\frac{C_{10}\tau_{B(y, d(x,y) - \wt u_n)}}{\phi(x,2\wt u_n)} \Big) \bigg] \bigg)^{\lfloor K\log n \rfloor -1}  \exp \Big( \frac{C_{10}\delta(1-\delta)^{n+1}}{\phi(x,2 \wt u_n)} \Big) \nn\\
	&\le e^6\exp \Big( -3K \log n +   \delta \kappa_{\delta,m+1}  \frac{C_{10} \phi(x,\wt u_n) \log |\log \phi(x,\wt u_n)|}{\phi(x,2 \wt  u_n)} \Big)\nn\\
	&\le e^{6}\exp \Big(\big(-3 K+2C_{10} \delta \kappa_{\delta,m+1} \big)\log  n\Big).
\end{align}
Besides, for all $n$ large enough, using the strong Markov property, \eqref{e:large_deviation_1}, the definition \eqref{e:def_u_n}, \eqref{e:loglog_u_n}, \eqref{e:scaling-zero} and the fact that $u_n \ge 4\wt u_n$ follows from \eqref{e:wtu_n_u_n}, we get that 
\begin{align}\label{4.1.2-2}
	& \P^x\left(\wt O_n \setminus O_n : \wt \tau^{n,2}_{\lfloor K\log n \rfloor} < (1-\delta)^n<\tau_{B(x,r_0)}\right)\le \P^x\left(\,\wt \tau^{n,2}_{\lfloor K\log n \rfloor} < (1-\delta)^n\right)  \nn\\
		&\le   \P^x 
		\bigg(\sum_{\ell=1}^{\lfloor K\log n \rfloor-1}\tau_{B(X_0, u_n - 2\wt u_n)} \circ \theta_{\wt \sigma^n_{b_\ell + 1}}  \le  \wt \tau^{n,2}_{\lfloor K\log n \rfloor} < (1-\delta)^n \bigg) \nn\\
	&\le \P^x 
	\bigg(\sum_{\ell=1}^{\lfloor K\log n \rfloor-1}\tau_{B(X_0, u_n - 2\wt u_n)} \circ \theta_{\wt \sigma^n_{b_\ell + 1}}  < (1-\delta)^n \bigg) \nn \\
	&\le \exp \Big( \frac{C_{10}  (1-\delta)^{n}}{\phi(x,2(u_n - 2\wt u_n))} \Big) \E^x\bigg[\exp \Big(-\frac{C_{10}}{\phi(x,2(u_n - 2\wt u_n))} \sum_{j=1}^{\lfloor K\log n \rfloor-1}\tau_{B(X_0, u_n - 2\wt u_n)} \circ \theta_{\wt \sigma^n_{b_\ell +1}}\Big) \bigg] \nn \\
	&\le \bigg(\sup_{y \in B(x, \wt u_n)}  \E^y\bigg[\exp \Big(-\frac{C_{10}\tau_{B(y, u_n-2\wt u_n)}}{\phi(x,2(u_n-2\wt u_n))} \Big) \Big] \bigg)^{\lfloor K\log n \rfloor -1}   \exp \Big( \frac{C_{10}  (1-\delta)^{n}}{\phi(x,2(u_n - 2\wt u_n))} \Big) \nn\\
	&\le e^6 \exp \bigg(-3K \log n +  \kappa_{\delta, m}   \frac{C_{10} \phi(x, u_n) \log |\log \phi(x,  u_n)|}{\phi(x,  u_n)}  \bigg)   \le e^6 \exp \Big(\big( -3K + 2C_{10} \kappa_{\delta, m} \big) \log n \Big). 
\end{align} 

Finally, using \eqref{4.1.1}, \eqref{4.1.2},  \eqref{4.1.2-1},  \eqref{4.1.2-2} and Proposition \ref{p:EP}(i), we deduce that for all $n$ large enough,
\begin{align}\label{4.1.5}
	\P^x(\wt O_n \setminus O_n) &\le 	\P^x\Big(\wt O_n\setminus O_n  :  \wt \tau^n_{\,2\lfloor K\log n \rfloor}\ge (1-\delta)^n\Big)  + 	\P^x\Big(\wt O_n\setminus O_n  :  \wt \tau^n_{\,2\lfloor K\log n \rfloor}< (1-\delta)^n < \tau_{B(x,r_0)}\Big)  \nn\\
	& \quad + 	\P^x\big( \tau_{B(x,r_0)} \le (1-\delta)^n\big)  \nn\\
	& \le\exp \Big( \big(2C_9 K - c_1 L \delta  \kappa_{\delta,m+1}\big)\log n\Big) +  e^{6}\exp \Big(\big(-3 K+2C_{10} \delta \kappa_{\delta,m+1} \big)\log  n\Big)  \nn\\
	& \quad  + e^6 \exp \Big(\big( -3K + 2C_{10} \kappa_{\delta, m} \big) \log n \Big) +  c_2 \phi(x,r_0)^{-q} (1-\delta)^{nq}.
\end{align}
Note that  the constant $L$ defined as \eqref{e:def_L} is  equal to $(2C_9+1)C_{10}/c_1$. Take 
$$\kappa_{\delta,m+1}= 2\delta^{-1}C_{10}^{-1} \vee \delta^{-1} \kappa_{\delta,m} \vee 4^{1+\beta_2}  C_U  \kappa_{\delta,m} \qquad \text{and} \qquad K = C_{10}\delta \kappa_{\delta,m+1}. $$ Then one can check that
\begin{equation}\label{e:induction_last} 
	\begin{cases}
		2C_9K - c_1L\delta \kappa_{\delta,m+1} = (2C_9C_{10}-c_1L)\delta  \kappa_{\delta,m+1} = -C_{10} \delta \kappa_{\delta,m+1}  \le -2, \\[3pt]
		-3K+ 2C_{10} \delta \kappa_{\delta,m+1} = -C_{10} \delta \kappa_{\delta,m+1} \le -2, \\[3pt]
		-3K + 2C_{10}\kappa_{\delta, m} \le  - 3K + 2C_{10}\delta \kappa_{\delta,m+1}\le -2. 
	\end{cases} 
\end{equation}
It now follows from \eqref{4.1.5} that $\P^x(\wt O_n \setminus O_n)$ is summable. The proof is complete. \qed

\noindent	\textbf{Proof of Theorem \ref{t:generalocc-SP}.} The result follows from  Lemmas \ref{l:occ-lem-1}(i), \ref{l:occ_lower}(i), \ref{l:generalocc-SP-4} and \ref{l:generalocc-SP-3}. \qed

\medskip

Theorem \ref{t:generalocc2} can be proved by a   similar way to that of Theorem \ref{t:generalocc-SP} with some modifications. 

Recall that  $\vp$ is the function defined in \eqref{e:molli1}.
Before giving the proof of Theorem \ref{t:generalocc2}, we first give a lemma about regularity of the function $\vp$.

\begin{lem}\label{l:phiinf}
	 If \ref{B2}  holds, then  for any $a>0$,
$
		\lim_{r \to \infty} {\vp(r + a)}/{\vp(r)} = 1.
$
\end{lem}
\pf Let $a>0$. Obviously, $	\liminf_{r \to \infty} \vp(r + a)/\vp(r) \ge 1$. On the other hand, by the mean value theorem,  \eqref{e:molli2} and \eqref{e:scaling-infty}, we obtain
\begin{align*}
	&\limsup_{r \to \infty}\frac{\vp(r+a) - \vp(r)}{\vp(r)} \le  c_1\limsup_{r \to \infty} \frac{\sup_{s \in (r, r+ a)} as^{-1}\phi(s)}{\phi(r)} \le c_2\limsup_{r \to \infty}  \frac{a}{r} \bigg( \frac{r+a}{r} \bigg)^{\beta_2}  = 0.
\end{align*}
Hence, we get $\limsup_{r \to \infty} \vp(r + a)/\vp(r) \le 1$ and  finish the proof.
\qed

\medskip 

\noindent	\textbf{Proof of Theorem \ref{t:generalocc2}.}  For $x \in M$, $\kappa>0$ and $a \in (0,1]$, define
\begin{equation*}
	E_{x,\kappa}(a)=\bigg\{\omega \in \Omega: \limsup_{r \to \infty}\frac{\sU\big(B(x,r), \kappa \vp(r) \log \log \vp(r) \big)(\omega)}{\kappa\vp(r)\log\log \vp(r) } <  a\bigg\}
\end{equation*}
and $g(x,\kappa)=\inf\{a \in (0,1]:\P^x(E_{x,\kappa}(a))=1\}$. Since $\lim_{r \to \infty} \vp(r) = \infty$ by \eqref{e:molli2} and Remark \ref{r:ass-B}(i),  one can see that  $E_{x,\kappa}(a)$ is shift-invariant for all $a \in (0,1]$. Thus, by Proposition \ref{p:law01}, for each fixed $a \in (0,1]$, it holds that either $\P^z(E_{x, \kappa }(a)) = 1$ for all $z \in M$, or $\P^z(E_{x,\kappa}(a)) = 0$ for all $z \in M$.  Therefore, for every $x \in M$ and $\kappa>0$,
\begin{equation}\label{e:g(x,k)}
	\limsup_{r \to \infty} \frac{\sU\big(B(x, r),\, \kappa\vp(r) \log \log \vp(r)\big)}{\kappa \vp (r) \log \log \vp(r)} = g(x,\kappa),~~\quad
	\P^z\mbox{-a.s.}, ~~\forall z \in M.
\end{equation}
 By the monotone property of occupation times and  Lemma \ref{l:phiinf}, it holds that for any $x,y \in M$,
\begin{align}\label{e:compareocc}
	g(x,\kappa) &\le \limsup_{r \to \infty} \frac{\sU\big(B(y, r+d(x,y)), \,\kappa\vp( r+d(x,y)) \log \log \vp(r+d(x,y)) \big)}{\kappa\vp(r) \log \log \vp(r) } \nn\\
	& \le g(y,\kappa) \limsup_{r \to \infty} \frac{\vp(r+d(x,y))\log \log \vp (r+d(x,y))}{\vp(r) \log \log \vp(r) } = g(y,\kappa).
\end{align}
Hence,  $g(x,\kappa)$ is equal to $g(y,\kappa)$ by symmetry so that is independent of $x$. Define $f_\infty(\kappa)=g(o,\kappa)$. By  \eqref{e:g(x,k)}, $f_\infty$ satisfies \eqref{e:generalocc2_1}. Moreover, by Lemmas \ref{l:occ-lem-1}(ii) and \ref{l:occ_lower}(ii),  $f_\infty$ is non-increasing and  there exists $\kappa_\infty>0$ such that $f_\infty(\kappa)=1$ for all $\kappa \le \kappa_\infty$. Thus it remains to show $\lim_{\kappa \to \infty}f_\infty(\kappa)=0$. To do this,  we  follow the proofs of Lemmas \ref{l:generalocc-SP-4},  \ref{l:generalocc-SP-2} and \ref{l:generalocc-SP-3}.

 For $\kappa>0$ and  $\delta \in (0,1/8)$,  let $(w^{\kappa, \delta }_n)_{n \ge N}$ be an increasing sequence such that (cf. \eqref{e:def_u_n})
\begin{equation}\label{e:def_w_n}
	\kappa \vp(w^{\kappa,\delta }_n)\log \log \vp(w^{\kappa, \delta}_n) = (1+\delta)^n \quad \text{for all} \; n\ge N.
\end{equation}
By assuming that $N$ is large enough, similar to  \eqref{e:loglog_u_n}, we get from \eqref{e:def_w_n} that 
\begin{equation}\label{e:loglog_w_n}
	2^{-1}\log n \le \log \log \vp(w^{\kappa, \delta}_n) \le 2\log n  \quad \text{for all} \;\; n\ge N.
\end{equation}
For $\kappa>0$, $\delta \in (0, 1/8)$  and $m \in \N$, we set
$$
\sO_n(\kappa, \delta,m):=\left\{\, \sU\big(B(o,  w^{\kappa, \delta}_n),  (1+\delta)^n \big) \ge (1+\delta)^{n-1}(1-m\delta)\right\}, \quad n \ge N.$$
Let
\begin{equation*}
	L':= 4^{1+\beta_2} C_U' (2C_9+1)C_{10}/C_8,
\end{equation*}
where $\beta_2,C_U'$ are the constants from \eqref{e:scaling-infty*} and  $C_8, C_9, C_{10}$ from  Lemma \ref{l:mean-exponential}(ii).  Cf. \eqref{e:def_L}. 
We will prove in Lemma \ref{l:generalocc-SP-5} below that there exists a constant $\delta_0 \in (0,1/8)$ such that for any $\delta \in (0, \delta_0]$ and $m \in \N$ satisfying $(m+L') \delta \le 1$, there exists a constant $\kappa_{\delta,m}>0$  such that 
\begin{equation}\label{e:sOni}
	\sum_{n=1}^\infty \P^x (\sO_n(\kappa_{\delta, m},\delta,m))	<\infty.
\end{equation}
Assume \eqref{e:sOni} for the moment. Choose  any  $\delta\in (0, \delta_0/(L'+1))$, and set $m:=\lfloor \delta^{-1}-L'\rfloor$ and $\lambda:=\kappa_{\delta,m}$. Then, similar to \eqref{e:generalocc_last}, using  the Borel-Cantelli lemma, we get from   \eqref{e:sOni}  that $\P^o$-a.s.,
\begin{align*}
\limsup_{r \to \infty} \frac{\sU(B(o, r), \lambda\vp(r)\log \log \vp(r))}{\lambda\vp(r)\log \log \vp(r)} 	& \le \limsup_{n \to \infty} \frac{\sU\big(B(o,  w_{n}^{\lambda, \delta}), \lambda\vp(w_{n}^{\lambda, \delta})\log \log \vp (w_{n}^{\lambda, \delta})\big)}{\lambda\vp(w_{n-1}^{\lambda, \delta})\log\log \vp(w_{n-1}^{\lambda, \delta})}\\
	&=\limsup_{n \to \infty} \frac{\sU\big(B(o,  w_{n}^{\lambda, \delta}),  (1+\delta)^{n}\big)}{(1+\delta)^{n-1}} \le  1- m \delta \le (L'+1)\delta.
\end{align*}
Thus, $f_\infty(\lambda) \le (L'+1)\delta$. By the monotone property of $f_\infty$, since $\delta$ can be chosen arbitrarily small, we finish  the proof.\qed

\begin{lem}\label{l:generalocc-SP-5}
	Under the setting of Theorem \ref{t:generalocc2},   \eqref{e:sOni} is true.
\end{lem}
\pf Let $\delta_0 \in (0,1/8)$ and $\lambda_1>0$ be constants chosen later. Choose any $\delta \in (0, \delta_0]$, and write $w_n=w_n^{\lambda_1,\delta}$ and $\sO_n=\sO_n(\lambda_1,\delta,1)$.  Then we define
\begin{equation*}
	\sigma_1^n:=0, \;\; 	\tau^n_{j}:= \inf\{s >\sigma^n_{j}: X_s \notin  B(o, 2w_n)  \} \quad \text{and} \quad \sigma^n_{j+1}:= \inf\{s >\tau^n_{j}: X_{s} \in B(o,  w_n)  \}.
\end{equation*}
 By following the calculations in \eqref{e:On-tau-le-log} and \eqref{e:generalocc-eq1}, we get from Lemma \ref{l:mean-exponential}(ii) that for all $n$ large enough,
\begin{align*}
	&\P^o\left(\sO_n  :  \tau^n_{\lfloor \log n \rfloor} \ge (1+\delta)^n\right) 	\le \P^o 
	\bigg(\sum_{j=1}^{\lfloor \log n \rfloor}\tau_{B(o, 2 w_n)} \circ \theta_{\sigma^n_j}  \ge  (1+\delta)^{n-1}(1-\delta) \bigg)\nn\\
	&\le  \exp \Big(- \frac{C_8(1+\delta)^{n-2}(1-\delta)}{\vp(4  w_n)} \Big)\E^o\bigg[\exp \Big(\frac{C_8}{\vp(4w_n)} \sum_{j=1}^{\lfloor \log n \rfloor}\tau_{B(o, 2 w_n)} \circ \theta_{\sigma^n_j}\Big) \bigg]\\
	& \le  \exp \Big(- \frac{C_8(1+\delta)^{n-2}(1-\delta)}{\vp(4  w_n)} \Big)  \bigg(\sup_{y \in  B(o, 2 w_n)} \E^y\bigg[\exp \Big(\frac{C_8\tau_{B(y, 4w_n)}}{\vp(4w_n)} \Big) \bigg] \bigg)^{\lfloor \log n \rfloor }\\
	&\le  \exp \Big( C_9 \log n - \frac{C_8(1+\delta)^{n-2}(1-\delta)}{\vp(4  w_n)} \Big).
\end{align*}
Note that  $(1+\delta)^{-2}(1-\delta)>1/2$ since $\delta<1/8$. By the definition \eqref{e:def_w_n} and \eqref{e:scaling-infty*}, it follows that for all $n$ large enough, 
\begin{align}\label{e:sO_n_largetau}
	\P^o\Big(\sO_n  :  \tau^n_{\lfloor \log n \rfloor} \ge (1+\delta)^n\Big) &\le  \exp \Big(  C_9 \log n - \lambda_1\frac{C_8 \vp(w_n) \log \log \vp(w_n) }{2 \vp (4w_n)}\Big) \nn\\
	& \le  \exp \Big( \big(C_9  -c_1 \lambda_1 \big) \log n\Big),
\end{align}
where $c_1:=4^{-1-\beta_2} C_8/C_U'$. Next, by similar calculations to that for obtaining  \eqref{e:generalocc-eq2}-\eqref{e:On_large_Tn}, we get that for  $p:=\up^{-1/3}>1$ and  all $n$ large enough, 
\begin{align}\label{e:defc_2}
	&\P^o\Big(\sO_n  :  \tau^n_{\lfloor \log n \rfloor} < (1+\delta)^n<\tau_{B(o, w_n^p)}\Big) \nn\\
	&\le \P^x 
	\bigg(\sum_{j=1}^{\lfloor \log n \rfloor-1}\tau_{B(X_0, w_n)} \circ \theta_{\tau^n_j}  \le 2\delta(1+\delta)^{n-1}, \;  X({\tau_j^n}) \in B(o,w_n^{p}),\, 1\le j \le \lfloor \log n \rfloor-1 \bigg)\nn\\
	&\le  \exp \Big( \frac{2C_{10}\delta (1+\delta)^{n-1}}{\vp(2 w_n)} \Big) \, \E^o\bigg[\exp \Big(-\frac{C_{10}}{\vp(2w_n)} \sum_{j=1}^{\lfloor \log n \rfloor-1}\tau_{B(X_0, w_n)} \circ \theta_{\tau^n_j}\Big) \nn\\
	&\hspace{55mm}:    X({\tau_j^n}) \in B(o,w_n^{p}),\, 1\le j \le \lfloor \log n \rfloor-1\bigg] \nn\\
	& \le  \exp \Big( \frac{2C_{10}\delta (1+\delta)^{n-1}}{\vp(2  w_n)} \Big)\bigg(\sup_{y \in  B(o, w_n^{p})}  \E^y\bigg[\exp \Big(-\frac{C_{10}\tau_{B(y, w_n)}}{\vp(2w_n)} \Big) \bigg] \bigg)^{\lfloor \log n \rfloor -1}\nn\\
	&\le e^6 \exp \Big( -3\log n +  \delta \lambda_1  \frac{2C_{10} \vp(w_n) \log \log \vp(w_n)}{(1+\delta)\vp(2 w_n)}   \Big)\le e^6\exp \Big(\big(-3 +4C_{10} \delta\lambda_1  \big)\log  n\Big).
\end{align}
Lastly, using Proposition \ref{p:EP}(ii) (with $\up_1= \up p$), \eqref{e:loglog_w_n} and \eqref{e:scaling-infty*}, we get that for all $n$ large enough,
\begin{align}\label{e:out-wp}
	\P^o\big(\tau_{B(o,w_n^p)} \le (1+\delta)^n \big) &= 	\P^o\big(\tau_{B(o,w_n^p)} \le \lambda_1 \vp(w_n) \log \log \vp(w_n) \big) \nn\\
	&\le c_2 \lambda_1  \frac{\vp(w_n) \log \log \vp(w_n) }{\vp(w_n^p)}\le 2c_2\lambda_1C_L'^{-1} w_n^{-\beta_1(p-1)} \log n.
\end{align}
Using \eqref{e:loglog_w_n} and \eqref{e:scaling-infty*}, we see from the definition \eqref{e:def_w_n} that for all $n$ large enough,
\begin{equation}\label{e:out-wp-2}
	w_n \ge c_3\bigg( \frac{\vp(w_n)}{\vp(2R_\infty)}\bigg)^{1/\beta_2} \ge c_3 \bigg( \frac{(1+\delta)^n}{2 \lambda_1  \vp(2R_\infty)\log n}\bigg)^{1/\beta_2},
\end{equation}
which yields that $\sum_{n=16}^\infty w_n^{-\beta_1(p-1)} \log n<\infty$.   Take $\lambda_1=(2+C_9)/c_1$ and $\delta_0=1/(4C_{10}\lambda_1)$. Then we get from \eqref{e:sO_n_largetau}, \eqref{e:defc_2} and \eqref{e:out-wp} that for all $n$ large enough, 
$$\P^o(\sO_n) \le (1+e^6)n^{-2} +  2c_2\lambda_1C_L^{-1} w_n^{-\beta_1(p-1)} \log n$$ and hence $\sum_{n=1}^\infty \P^o(\sO_n)<\infty$. This proves the base step.

For the next step, let $m \ge 1$ be such that  $(m+1+L')\delta  \le 1$ and suppose that \eqref{e:sOni} holds with $m$.  Let $\kappa_{\delta,m+1} \ge 4^{1+\beta_2}C_U'\kappa_{\delta,m}$ be a constant  chosen later. Denote by
$$
v_n= w_n^{\kappa_{\delta, m},\,\delta}, \quad \wt v_n = w_n^{\kappa_{\delta, m+1},\,\delta}, \quad U_n=\sO_n(\kappa_{\delta, m},\delta,m)  \quad \text{and} \quad  \wt U_n=
\sO_n(\kappa_{\delta, m+1},\delta,m+1).
$$ 
As in the proof of Lemma \ref{l:generalocc-SP-3}, our task is to show  $\sum_{n=1}^\infty \P^o(\wt U_n \setminus U_n) < \infty$.  By similar calculations to \eqref{e:wtu_n_u_n}, we see that $4 \wt v_n \le v_n$  for all $n$ large enough. Let
\begin{align*}
	\wt \sigma_1^n:=0, \;\; \wt 	\tau^n_{j}:= \inf\{s >\wt \sigma^n_{j}: X_s \notin  B(o, 2\wt v_n)  \} \quad \text{and} \quad \wt  \sigma^n_{j+1}:= \inf\{s >\wt \tau^n_{j}: X_{s} \in B(o,  \wt v_n)  \}.
\end{align*}

By analogous arguments as the ones for obtaining  \eqref{4.1.1}, using \eqref{e:def_w_n}, \eqref{e:scaling-infty*} and the assumption that $(m+1+L')\delta \le 1$, we get that for any $K>0$ and all $n$ large enough,
\begin{align*}
	&\P^o\Big(\wt U_n \setminus U_n  :  \wt \tau^n_{2\lfloor K\log n \rfloor} \ge (1+\delta)^n\Big)\nn\\
	&\le \P^o
	\bigg(\sum_{j=1}^{2\lfloor K\log n \rfloor} \tau_{B(o, 2\wt v_n)} \circ \theta_{\wt \sigma^n_j}  \ge  (1+\delta)^{n-1} (1- \delta - m\delta)  \bigg)\nn\\
	&\le  \exp \Big(- \frac{C_8(1+\delta)^{n-1}(1-\delta -  m\delta)}{\vp(4\wt v_n)} \Big) \,\E^o\bigg[\exp \Big(\frac{C_8}{\vp(4 \wt v_n)} \sum_{j=1}^{2\lfloor K\log n \rfloor}\tau_{B(o, 2 \wt v_n)} \circ \theta_{\wt \sigma^n_j}\Big) \bigg]\\
	& \le \cdots \le  \exp \Big(- \frac{C_8(1+\delta)^{n-1}(1-\delta -  m\delta)}{\vp(4\wt v_n)} \Big) \bigg(\sup_{y \in  B(o, 2\wt v_n)} \E^y\bigg[\exp \Big(\frac{C_8 \tau_{B(y, 4\wt v_n)}}{\vp(4\wt v_n)} \Big) \bigg] \bigg)^{2\lfloor K\log n \rfloor }\\
	&\le  \exp \bigg( 2C_9 K \log n - \kappa_{\delta,m+1}  \frac{C_8(1+\delta)^{-1}(1-\delta -m\delta)  \vp(\wt v_n) \log \log \vp(\wt v_n)}{\vp(4\wt v_n)} \bigg)\\
	& \le \exp \Big( \big(2C_9 K - c_1 L' \delta  \kappa_{\delta,m+1}\big) \log n\Big),
\end{align*}
where $c_1>0$ is the constant from \eqref{e:sO_n_largetau}.

Next, with a slight abuse of notation, we define  $(a_\ell)_{\ell\ge 0}$ and $(b_\ell)_{\ell\ge 0}$ as \eqref{e:def-a_n} and \eqref{e:def-b_n} with $v_n$ and $\wt v_n$ instead of $u_n$ and $\wt u_n$ therein respectively.  
Then we write  $\wt\tau^{n,1}_\ell= \wt\tau^n_{a_\ell}$ and $\wt\tau^{n,2}_\ell = \wt \tau^n_{b_\ell}$. Similar to  \eqref{e:O_n_minus_O_n}, one can see that on the event $\wt U_n \setminus U_n$,  
 $\sU(B(y,\wt v_n), (1+\delta)^n) \le \delta(1+\delta)^{n-1}$  for all $y \in B(o,v_n - \wt v_n) \setminus B(o,2\wt v_n)$. Thus, using the strong Markov property and Lemma \ref{l:mean-exponential}(ii), we get  that for any $K>0$ and  all $n$ large enough,
\begin{align*}
	&\P^o\left(\wt U_n \setminus U_n:  \wt \tau^{n,1}_{\lfloor K \log n \rfloor} < (1+\delta)^n<\tau_{B(o, \wt v_n^p)}\right) \\
	&\le \P^o
	\bigg(\sum_{j=1}^{\lfloor K\log n \rfloor-1}\tau_{B(X_0, \wt v_n)} \circ \theta_{\wt \tau^{n,1}_j}  \le \delta(1+\delta)^{n-1}, \;  X_{\wt \tau_j^{n,1}} \in B(o,\wt v_n^{\,p}),\, 1\le j \le \lfloor K\log n \rfloor-1 \bigg)\\
	& \le   \exp \Big( \frac{C_{10}\delta(1+\delta)^{n-1}}{\vp(2 \wt v_n)} \Big) \, \bigg(\sup_{y \in  B(o, \wt v_n^{\,p})}  \E^y\bigg[\exp \Big(-\frac{C_{10}\tau_{B(y, \wt v_n)}}{\vp(2\wt v_n)} \Big) \bigg] \bigg)^{\lfloor K\log n \rfloor -1}  \\
	& \le e^6 \exp \Big( -3K \log n +  \delta \kappa_{\delta,m+1} \frac{C_{10} \vp( \wt v_n) \log \log \vp(\wt v_n)}{(1+\delta)\vp(2 \wt v_n)}\Big) \le e^6 \exp \Big(\big(-3 K+ 2C_{10} \delta \kappa_{\delta,m}  \big)\log  n\Big).
\end{align*}
Moreover, by a similar argument to the one for obtaining  \eqref{4.1.2-2}, we see that for all $n$ large enough,
\begin{align*}
	&\P^o\left(\wt U_n \setminus U_n : \wt \tau^{n,2}_{\lfloor K\log n \rfloor} < (1+\delta)^n<\tau_{B(o, \wt v_n^p)}\right)\\
	&\le   \P^o
	\bigg(\sum_{\ell=1}^{\lfloor K\log n \rfloor-1}\tau_{B(X_0, v_n - 2\wt v_n)} \circ \theta_{\wt \sigma^n_{b_\ell + 1}}  < (1+\delta)^n \bigg)  \\
		&\le \exp \Big( \frac{C_{10} (1+\delta)^{n}}{\vp(2(v_n - 2\wt v_n))} \Big) \E^o\bigg[\exp \Big(-\frac{C_{10}}{\vp(2(v_n - 2\wt v_n))} \sum_{j=1}^{\lfloor K\log n \rfloor-1}\tau_{B(X_0, v_n - 2\wt v_n)} \circ \theta_{\wt \sigma^n_{b_\ell +1}}\Big) \bigg] \nn \\
	&\le \bigg(\sup_{y \in B(o, \wt v_n)}  \E^y\bigg[\exp \Big(-\frac{C_{10}\tau_{B(y, v_n-2\wt v_n)}}{\vp(2 (v_n-2\wt v_n))} \Big) \bigg] \bigg)^{\lfloor K\log n \rfloor -1}   \exp \Big( \frac{C_{10}  (1+\delta)^{n}}{\vp(2(v_n - 2\wt v_n))} \Big) \nn\\
	&\le e^6 \exp \bigg( -3K \log n + \kappa_{\delta, m}   \frac{C_{10} \vp( v_n) \log \log \vp(  v_n)}{\vp(  v_n)}  \bigg)   \le e^6 \exp \Big(\big( -3K + 2C_{10} \kappa_{\delta, m} \big) \log n \Big). 
\end{align*}

Lastly, by similar calculations to the ones for obtaining  \eqref{e:out-wp} and \eqref{e:out-wp-2}, we get  that  for all $n$ large enough,
\begin{align*}
\P^o\big(\tau_{B(o, \wt v_n^p)} \le (1+\delta)^n \big) &\le c_4 \kappa_{\delta, m+1} \frac{\vp(\wt v_n) \log \log \vp (\wt v_n)}{\vp(\wt v_n^p)} \le 2c_4 \kappa_{\delta, m+1} C_L'^{-1} \wt v_n^{-\beta_1(p-1)} \log n\\
& \le c_5 (1+\delta)^{-n\beta_1(p-1)/\beta_2} (\log n)^{1-1/\beta_2}
\end{align*}
and hence $\sum_{n=N}^\infty \P^o\big(\tau_{B(o, \wt v_n^p)} \le (1+\delta)^n \big)  <\infty$.

Notice that $L'=(2C_9+1)C_{10}/c_1$.  Take $\kappa_{\delta,m+1}:= 2 \delta^{-1}C_{10}^{-1} \vee \delta^{-1}\kappa_{\delta,m} \vee 4^{1+\beta_2}  C_U' \kappa_{\delta,m}$ and $K =C_{10} \delta \kappa_{\delta,m+1}$. Then by the above four displays  and \eqref{e:induction_last}, we deduce that
$ \sum_{n=N}^\infty \P^o(\wt U_n \setminus U_n)< \infty$. The proof is complete. \qed

We end this section with  the proof of Theorem \ref{t:generalocc3}.

\medskip

\noindent	\textbf{Proof of Theorem \ref{t:generalocc3}.} For  $x \in M$ and $a>0$, we  define
$$E_x'(a):=\bigg\{ \limsup_{r \to \infty}\frac{\sU(B(x,r), \infty)}{\vp(r)\log\log \vp(r) } <  a   \bigg\}.$$
Since $X$ is conservative by Proposition \ref{p:EP}(ii), $E_x'(a)$ is shift-invariant. Moreover,  using Lemma \ref{l:phiinf}, similar to \eqref{e:compareocc},  one can check that $E_x'(a)$ is independent of $x$.
 Thus, by Proposition \ref{p:law01},  it suffices to show that there are constants $p_2 \ge p_1>0$ such that
\begin{equation}\label{e:generalocc3}
	\limsup_{r \to \infty}\frac{\sU(B(o,r), \infty)}{\vp(r)\log\log \vp(r)} \in [p_1,p_2], \quad \P^o\text{-a.s.}
\end{equation}
The lower bound in \eqref{e:generalocc3} comes from \eqref{e:generalocc2_1} and the monotonicity of occupation times.

Write $\oB(r):=\{x \in M: d(o,x) \le r\}$.  Observe that for all $r \ge (2R_\infty)^{1/(1-\up)}$ and $x \in \oB(r)$,
$$
R_\infty \d(x)^\up \le R_\infty (r+1)^\up < 2 R_\infty r^\up  \le r.
$$
Hence, by \ref{NDU} and   \VDi, we see that  for all  $r \ge (2R_\infty)^{1/(1-\up)}$, $x,y \in \oB(r)$ and $t \ge \vp(r)$, 
\begin{equation}\label{e:UHKD}
	p(t,x,y) \le  \frac{c_1}{V(x, \vp^{-1}(t)) \wedge V(y, \vp^{-1}(t))} \le  \frac{c_2}{V(x, 2\vp^{-1}(t)) \wedge V(y, 2\vp^{-1}(t))} \le \frac{c_2}{V(o, \vp^{-1}(t))}.
\end{equation}
Using \eqref{e:UHKD}, \VDi \ and \eqref{e:scaling-infty*} twice, since $\beta_2<d_1$, we get that for all $r \ge (2eR_\infty)^{1/(1-\up)}$ and $x \in \overline B(r)$, 
\begin{align*}
	\E^x\big[\,\sU(\oB(r),\infty)\big] &=\int_0^{\vp(2r)} \int_{\oB(r)}p(t,x,y)\mu(dy) dt +  \int_{\vp(2r)}^\infty \int_{\oB(r)}p(t,x,y)\mu(dy) dt\\
	& \le \vp(2r) + c_2\int_{\vp(2r)}^\infty \frac{ V(o,2r)}{ V(o, \vp^{-1}(t))}dt  \le\vp(2r) + c_3 \int_{\vp(2r)}^\infty \bigg(\frac{2r}{\vp^{-1}(t)}\bigg)^{d_1} dt \\
	& \le \vp(2r) +  c_4\int_{\vp(2r)}^\infty\bigg(\frac{\vp(2r)}{t}\bigg)^{d_1/\beta_2}  dt  = c_5\vp(2r) \le c_6 \vp(r/e).
\end{align*}
Thus, by Lemma \ref{l:suprema},  we deduce that  for all $r \ge (2eR_\infty)^{1/(1-\up)}$, 
$$\sup_{x \in M}  \E^x\big[\,\sU( \oB(r), \infty)\big] = \sup_{x \in \oB(r)}  \E^x\big[\,\sU(\oB(r), \infty)\big]  \le  c_6 \vp(r/e).$$
By  Khasmin'skii's lemma (see \cite[Lemma 3]{Kh} or \cite[Lemma B.1.2]{Sb}), it follows  that
\begin{equation*}\sup_{x \in M}	\E^x \bigg[ \exp\bigg( \frac{\sU(\oB(e^{n+1}),\infty)}{2c_6\vp(e^{n})}\bigg)\bigg] \le 2 \quad \text{for all $n$ large enough}.
\end{equation*}
The lower inequality in \eqref{e:scaling-infty*} implies that $\log \log \vp(e^n) \ge 2^{-1} \log n$ for all $n$ large enough. Hence, by Markov inequality and the monotonicity of occupation times, we get that for all $n$ large enough,
\begin{align*}
	\P^o \big(\, \sU(B(o,e^{n+1}),\infty)> 8c_6\vp(e^n) \log\log \vp(e^n)\big) &\le \P^o \bigg(\,  \frac{\sU(B(o,e^{n+1}),\infty)}{4c_6\vp(e^n)}  \ge \log n\bigg)\\
	&\le n^{-2}\E^o \bigg[ \exp\Big( \frac{\sU(\oB(e^{n+1}),\infty)}{2c_6\vp(e^n)}\Big)\bigg] \le 2 n^{-2}.
\end{align*}
Finally, using the Borel-Cantelli lemma, we arrive at
\begin{align*}
	\limsup_{r \to \infty}\frac{\sU(B(o,r), \infty)}{\vp(r)\log\log \vp(r)} &\le \limsup_{n \to \infty} \sup_{e^n\le r\le e^{n+1}}\frac{\sU(B(o,r), \infty)}{\vp(r)\log\log \vp(r)} \\
	&\le  \limsup_{n \to \infty} \frac{\sU(B(o,e^{n+1}), \infty)}{\vp(e^n)\log\log \vp(e^n)} \le 8c_6, \qquad \P^o\text{-a.s.}
\end{align*}
which completes the proof. \qed

\section{Proof of Theorem  \ref{t:occ1}}\label{s:occ2}

In this section except in the proof of Corollary \ref{c:occ1}, we always assume that \eqref{e:recur} holds for some $x_0 \in M$, which implies that \eqref{e:recur}  holds with $x_0=o$ by Remark \ref{r:NDLrecur}.  For simplicity, let us denote $\lVert F \rVert_p$ for $\lVert F \rVert_{L^p(M,\mu)}$, $p\in\{1,\infty\}$.

\begin{lem}\label{l:occ1-1}
Let $a>0$.	For every $F \in L^1(M,\mu)$ with $\lVert F \rVert_1 \neq 0$,  the event 
$$E_F(a):=\bigg\{ \omega \in \Omega\,; \, \limsup_{t \to \infty} \frac{\sU(F,t)(\omega)/\lVert F\rVert_1}
{\Theta(t/\log\log \Theta(t)) \log\log\Theta(t)} <a \bigg\}$$
 is  shift-invariant.
\end{lem}
\pf   Define $g(t)=\Theta(t/\log\log \Theta(t)) \log\log\Theta(t)$. Then $\lim_{t \to \infty}g(t)=\infty$ by \eqref{e:recur}. Hence, $E_F(a)$ is a tail event. Moreover, we have that, for all $s > 0$ and $\omega \in \Omega$,
\begin{align*}
&\limsup_{t \to \infty} \frac{\sU(F,t) (\theta_s \omega) }{g(t)} 	=\limsup_{t \to \infty} \frac{\sU(F,t+s) - \sU(F,s) }{g(t)}(\omega) \\
&\le  \limsup_{t \to \infty} \frac{\sU(F,t)(\omega) + \sU(F,s)(\theta_t \omega)  }{g(t)} \le   \limsup_{t \to \infty} \frac{\sU(F,t)(\omega) + s \lVert F \rVert_\infty }{g(t)} = \limsup_{t \to \infty} \frac{\sU(F,t)( \omega) }{g(t)}
\end{align*}
and
\begin{align*}
	&\limsup_{t \to \infty} \frac{\sU(F,t+s) - \sU(F,s) }{g(t)}(\omega) \\
	& \ge  \limsup_{t \to \infty} \frac{\sU(F,t)- \sU(F,s) }{g(t)}(\omega)\ge   \limsup_{t \to \infty} \frac{\sU(F,t)(\omega) - s \lVert F \rVert_\infty }{g(t)} = \limsup_{t \to \infty} \frac{\sU(F,t) (\omega)}{g(t)}.
\end{align*}
Therefore,  $\limsup_{t \to \infty} \sU(F,t) (\theta_s \omega)/g(t) =\limsup_{t \to \infty} \sU(F,t)(\omega)/g(t)$.  \qed

\noindent \textbf{Proof of Theorem \ref{t:occ1}.} 
In Lemmas \ref{l:occ1-2} and \ref{l:occ1-3} below, we will prove  that there are constants $a_2 \ge a_1>0$ such that for every $ F \in \fM$ satisfying \eqref{e:F} and $\lVert F \rVert_1 \neq 0$,
\begin{equation}\label{e:occ_generalform_1}
	\limsup_{t \to \infty} \frac{\sU(F,t)/ \lVert F \rVert_1}{\Theta( t/\log \log \Theta(t)) \log \log  \Theta(t)} \in [a_1,a_2],\quad\,
	\P^o\mbox{-a.s.}
\end{equation}
Then using Proposition \ref{p:law01} and Lemma \ref{l:occ1-1}, we deduce \eqref{e:occ_generalform} from \eqref{e:occ_generalform_1}. \qed

The following lemma will be used in the proof for the upper bound in \eqref{e:occ_generalform_1}.
\begin{lem}\label{l:occ1-key}
	Suppose that \VDi,  \ref{B2} and \ref{NDU} hold. Then there exists a constant $K_0>0$ such that for every  $ F \in \fM$ satisfying \eqref{e:F} and $\lVert F \rVert_1 \neq 0$,  
\begin{equation}\label{e:meanocc}
	\sup_{w \in M}\E^w \big[\, \sU(F, t) \big] \le  K_0\lVert F \rVert_{1}\Theta(t) \quad  \text{for all} \;\; t>  T_0
\end{equation}
with some constant $T_0>0$.
\end{lem}
\pf  When $F \in \fMi$, let $r_1 > R_\infty^{1/(1-\up)}$ be a constant satisfying $\text{supp}[F] \subset B(o, r_1)$.
When $\beta_1>d_2$ and $F \in \fMg$ for $\gamma \in (\beta_1 d_2/(\beta_1 - d_2), \infty)$,  let $\delta_1 \in (0,1)$ be  a  constant such that 
\begin{equation}\label{e:def_delta1}
	d_2 \left( \frac{1}{\beta_1} + \frac{1}{\gamma} \right) + \beta_2\delta_1<1.
\end{equation}
Recall that $\overline B(r):=\{x \in M: d(o,x) \le r\}$.  We define for $t>0$,
\begin{equation*}
	\mathbf r(t)=\begin{cases}
		t^{\delta_1}\vp^{-1}(t)^{d_2/\gamma}, &\mbox{if} \;\; \gamma<\infty,\\
		r_1, &\mbox{if} \;\; \gamma=\infty,
	\end{cases} \;\;\quad F_{1,t}=F \cdot \1_{\overline B(\mathbf r(t))} \quad \text{and} \quad F_{2,t}=F-F_{1,t}.
\end{equation*} 
Observe that  for each $t>0$,
\begin{equation}\label{e:meanocc_decomp}
	\sup_{w \in M}\E^w \big[\, \sU(F, t) \big] \le \sup_{w \in M}\E^w \big[\, \sU(F_{1,t}, t) \big] + \sup_{w \in M}\E^w \big[\, \sU(F_{2,t}, t) \big].
\end{equation}

By \VDi, for all large $t$ such that $\mathbf r(t)^{1-\up}>R_\infty $, any $y \in \overline B(\mathbf r(t))$ and  $s > \vp(\mathbf r(t))$,  we have
$$ V(o,\vp^{-1}(s)) \le V(y, 2\vp^{-1}(s)) \le c_1 V(y, \vp^{-1}(s)).$$
Hence, by \ref{NDU}, it holds that  for all $t$ large enough and all $w \in  \overline B(\mathbf r(t))$, 
\begin{align*}
	\E^w \big[\, \sU(F_{1,t}, t) \big] &= \int_0^{\vp(\mathbf r(t))} \int_{\overline B(\mathbf r(t))} \, p(s, w, y)F(y) \mu(dy) ds + \int_{\vp(\mathbf r(t))}^t \int_{\overline B(\mathbf r(t))} \, p(s, w, y)F(y) \mu(dy) ds \\
	& \le \vp(\mathbf r(t))\lVert F \rVert_{\infty}+ c_2\int_{\vp(\mathbf r(t))}^t \int_{\overline B(\mathbf r(t))}  \left(\frac{1}{V(w, \vp^{-1}(s))} \vee \frac{1}{V(y, \vp^{-1}(s))} \right) F(y) \mu(dy)ds \\
	& \le \vp(\mathbf r(t))\lVert F \rVert_{\infty} + c_1c_2 \lVert F \rVert_{1}\int_{\vp(\mathbf r(t))}^t   \frac{ds}{V(o, \vp^{-1}(s))} \le \vp(\mathbf r(t))\lVert F \rVert_{\infty} + c_1c_2  \lVert F \rVert_{1}\Theta(t).
\end{align*}
Thus, by Lemma \ref{l:suprema}, it holds that for all $t$ large enough, 
\begin{equation}\label{e:meanocc_decomp1}
	\sup_{w \in M}\E^w \big[\, \sU(F_{1,t}, t) \big]=\sup_{w \in \overline B(\mathbf r(t))}\E^w \big[\, \sU(F_{1,t}, t) \big] \le  \vp(\mathbf r(t))\lVert F \rVert_{\infty} + c_1c_2 \lVert F \rVert_{1}\Theta(t).
\end{equation}

On the other hand, we have that, if $\gamma<\infty$, then 
\begin{align*}
	\sup_{w \in M}\E^w \big[\, \sU(F_{2,t}, t) \big] &=\sup_{w \in M} \int_0^{t} \int_{M \setminus \overline B(\mathbf r(t))} \, p(s, w, y)F(y) \mu(dy) ds \\
	&\le  (1+\mathbf r(t))^{-\gamma}\lVert \d(y)^\gamma F(y)\rVert_{\infty} \int_0^{t}ds  \sup_{w \in M}  \int_{M \setminus \overline B(\mathbf r(t))} \, p(s, w, y) \mu(dy) \\
	&\le  t(1+\mathbf r(t))^{-\gamma}\lVert \d(y)^\gamma F(y)\rVert_{\infty}
\end{align*}
and if $\gamma=\infty$, then $F_{2,t}=0$ so that $\sup_{w \in M}\E^w \big[\, \sU(F_{2,t}, t) \big] =0$.  Combining with  \eqref{e:meanocc_decomp1}, we get from  \eqref{e:meanocc_decomp}  that for all $t$ large enough,
$$\sup_{w \in M}\E^w \big[\, \sU(F, t) \big] \le c_1c_2  \lVert F \rVert_{1}\Theta(t)+  \vp(\mathbf r(t))\lVert F \rVert_{\infty} +\1_{\{\gamma<\infty\}} t(1+\mathbf r(t))^{-\gamma}\lVert \d(y)^\gamma F(y)\rVert_{\infty}.$$
Therefore, to obtain \eqref{e:meanocc}, it suffices to show that 
\begin{equation}\label{e:meanocc_dominant}
	\lim_{t \to \infty} \big(\vp(\mathbf r(t))+ \1_{\{\gamma<\infty\}} t(1+\mathbf r(t))^{-\gamma}\big)/\Theta(t)=0.
\end{equation}
If $\gamma=\infty$, then \eqref{e:meanocc_dominant} immediately follows since $\lim_{t \to \infty} \Theta(t)=\infty$. Now suppose that  $\gamma<\infty$ and $\beta_1>d_2$. By Lemma \ref{l:theta}(i) and \VDi, we have
\begin{align*}
	\lim_{t\to \infty} \frac{t(1+\mathbf r(t))^{-\gamma}}{\Theta(t)} &\le 2V(o,2R_\infty) \lim_{t\to \infty} \bigg(t^{-\delta_1\gamma} \vp^{-1}(t)^{-d_2}\frac{V(o,\vp^{-1}(t))}{V(o,2R_\infty)}\bigg)\\
	& \le c_3 R_\infty^{-d_2}V(o,2R_\infty)\lim_{t\to \infty}t^{-\delta_1\gamma}=0.
\end{align*}
Further, using Lemma \ref{l:theta}(i), \VDi \ and  \eqref{e:scaling-infty*} several times, we get from \eqref{e:def_delta1} that
\begin{align*}
	\lim_{t \to \infty}	\frac{\vp(\mathbf r(t))}{\Theta(t)} &\le 2\lim_{t \to \infty} \frac{ \vp\big(t^{\delta_1} \vp^{-1}(t)^{d_2/\gamma}\big) V(o,\vp^{-1}(t))}{t} =2\lim_{u \to \infty} \frac{\vp\big(\vp(u)^{\delta_1} u^{d_2/\gamma}\big) V(o,u)}{\vp(u) } \\
	&= 2V(o,2R_\infty)\lim_{u \to \infty} \bigg( \frac{\vp\big( \vp(2R_\infty)^{\delta_1} \vp(u)^{\delta_1} u^{d_2/\gamma} / \vp(2R_\infty)^{\delta_1}\big)}{\vp(u)} \frac{ V(o,u)}{V(o,2R_\infty)} \bigg)\\
	&\le c_4R_\infty^{-d_2}\lim_{u \to \infty}\frac{\vp\big( c_5u^{\beta_2\delta_1+d_2/\gamma} \big)}{\vp(u)} u^{d_2}\le c_6 \lim_{u \to \infty} u^{d_2- \beta_1(1-\beta_2\delta_1-d_2/\gamma)}=0.
\end{align*}
Thus, \eqref{e:meanocc_dominant} holds and we finish the proof. \qed

\begin{lem}\label{l:occ1-2}
	Under the setting of Theorem \ref{t:occ1}, the upper bound in  \eqref{e:occ_generalform_1} holds.
\end{lem}
\pf  By    Markov inequality and Lemma \ref{l:occ1-key}, it holds that for all $t \ge T_0$,
\begin{equation*}
	\sup_{w \in M} \P^w \Big( \sU(F,t)/\lVert F \rVert_{1} \ge 2K_0 \Theta(t) \Big) \le \frac{1}{ 2K_0\lVert F \rVert_{1} \Theta(t)} \sup_{w \in M}\E[\, \sU(F,t)]\le\frac{1}{2}.
\end{equation*}
Hence, using \eqref{e:recur} (with $x_0=o$) and Lemma \ref{l:theta}(ii),  we see that with $F(t)=\sU(F,t)/\lVert F \rVert_{1}$ and $g(t)=\Theta(t)$, conditions (1)--(3) of Lemma \ref{l:limsup} are satisfied. Now the upper bound in \eqref{e:occ_generalform_1} follows from Lemma \ref{l:limsup}. \qed

\begin{lem}\label{l:occ1-3}
	Under the setting of Theorem \ref{t:occ1},	the lower bound in \eqref{e:occ_generalform_1} holds.
\end{lem}
\pf	By \eqref{e:recur}, there exist  $N \ge 10$ and  an increasing sequence $(t_n)_{n \ge N}$ such that 
\begin{equation}\label{e:choosetn}
	\Theta(t_n/\log n) = \exp(n^2), \quad  n \ge N.
\end{equation} 
Set $u_n:= \log \log \Theta(t_n)$. Then  $u_n \ge  \log \log \exp(n^2) = 2 \log n$ for $n \ge N$. Further, by Lemma \ref{l:theta}(ii), we see that $u_n \le   \log \log \big(3(\log n)^{\log 3/\log 2} \exp(n^2)\big) \le \log \log \exp(n^3) = 3\log n$ for  $n \ge N$. Thus,
\begin{align}\label{e:usize}
	2\log n \le u_n \le 3\log n \quad \text{for all} \;\; n \ge N.
\end{align}	Using Lemma \ref{l:theta}(ii) twice, we get from  \eqref{e:choosetn} and \eqref{e:usize} that for all $n \ge N \ge 10$,
\begin{align}\label{e:gap}
	\Theta(n^{2}t_n) &\le 3(n^{2}\log n)^2 \Theta(t_n/\log n)=3(n^{2}\log n)^{2}e^{-2n-1} \Theta( t_{n+1}/\log (n+1)) \nn\\
	& \le 250(n^{2}\log n)^2 e^{-2n-1}\Theta(t_{n+1}/2u_{n+1}) \le  \frac{1}{2}	\Theta(t_{n+1}/2u_{n+1}).
\end{align}
In particular, since $\Theta$ is increasing, we see that $t_{n+1} \ge 100 t_n$ for all $n \ge N$.

Let  $r_F>0$ be  a constant such that $2\lVert  F \cdot \1_{\oB(r_F)}\rVert_{1} \ge \lVert F\rVert_{1}$ and set $F_0:= F \cdot \1_{\oB(r_F)}$.  By \ref{B4+}, there exist constants  $c_1>1$, $c_2>0$ such that for all $n$ large enough, all $w \in \oB( \vp^{-1}(nt_n))$ and $s> \vp(c_1 \vp^{-1}(nt_n))$,
\begin{equation}\label{e:NDL-F}
	\int_{\oB(r_F)}p^{B(o, c_1 \vp^{-1}(s))}(s,w,y) F_0(y)\mu(dy)\ge \frac{c_2 \lVert F_0 \rVert_1}{V(o,\vp^{-1}(s))}.
\end{equation}
Note that by \eqref{e:scaling-infty}, it holds that for all $n$ large enough,
\begin{equation*}
\vp(c_1 \vp^{-1}(nt_n)) \le c_2 n t_n \le n^2t_n.
\end{equation*}
Thus, using \eqref{e:NDL-F}, \eqref{e:gap}, Lemma \ref{l:theta}(ii) and the fact that $2\lVert F_0 \rVert_1 \ge \lVert F \rVert_1$, we get that  for all  $n$ large enough and any $w \in \oB(\vp^{-1}(nt_n))$,
\begin{align}\label{e:meanocc2}
	&\E^w \bigg[\,\sU\Big(F_0,\frac{t_{n+1}-t_n}{u_{n+1}}\Big) \bigg] \ge   \int_{\vp(c_1\vp^{-1}(nt_n))}^{t_{n+1}/(2u_{n+1})} \int_{\oB(r_F)} \,  p^{B(o, c_1\vp^{-1}(s))}(s,w,y)F_0(y)\mu(dy)  ds \nonumber\\
	&\ge c_2\lVert F_0\rVert_{1}\int_{n^2t_n}^{t_{n+1}/(2u_{n+1})}  \frac{ds }{V(o,\vp^{-1}(s))}  = c_2\lVert F_0 \rVert_{1}\big(\Theta( t_{n+1}/2u_{n+1} ) - \Theta(n^2t_n) \big) \nn\\[3pt]
	&\ge \frac{c_2}{2}\lVert F_0 \rVert_{1}\Theta( t_{n+1}/2u_{n+1} ) \ge \frac{c_2}{18} \lVert F_0 \rVert_{1}\Theta( t_{n+1}/u_{n+1} ) \ge \frac{c_2}{36} \lVert F \rVert_{1}\Theta( t_{n+1}/u_{n+1} ).
\end{align}
On the other hand, by Lemma \ref{l:occ1-key} and  Khasmin'skii's lemma (see \cite[Lemma 3]{Kh} or \cite[Lemma B.1.2]{Sb}), it holds that for all $n$ large enough,
\begin{equation*}
	\sup_{w \in M}\E^w \left[\,\sU\Big(F_0,\frac{t_{n+1}-t_n}{u_{n+1}}\Big)^2\right] \le  \sup_{w \in M}\E^w \left[\,\sU\Big(F,\frac{t_{n+1}}{u_{n+1}}\Big)^2\right]  \le \big(K_0\lVert F \rVert_{1} \Theta(t_{n+1}/u_{n+1})\big)^2.
\end{equation*}
Therefore, by combining the above with \eqref{e:meanocc2} and using the  Paley-Zygmund inequality, we deduce that  there are constants $c_3>0$ and $p \in (0,e^{-1})$ independent of $F$ such that for all $n$ large enough,
\begin{equation}\label{e:occlower}
	\inf_{w \in \oB(\vp^{-1}(nt_n) )} \P^w\bigg( \sU\Big(F_0, \frac{t_{n+1}-t_n}{u_{n+1}}\Big) \ge 2c_3 \lVert F \rVert_{1} \Theta(t_{n+1}/u_{n+1})\bigg) \ge p.
\end{equation}

 Set $\delta:=(3\log(1/p))^{-1}$ and define for $n \ge N$ and $i \ge 1$,
\begin{align*}
E_{n}&:= \bigg\{ \, \frac{\sU(F, t_{n})}{u_{n}\Theta(t_{n}/u_{n})  } \ge  \delta c_3  \lVert F \rVert_{1}  \bigg\},   \quad\;\;  	G_n:= \big\{ X_{t_n} \in \oB(\vp^{-1}(nt_n)) \big\}, \\[6pt]
H_{n,i}&:= \bigg\{ \exists s_1,...,s_i \in (0, (t_{n+1}-t_n)/u_{n+1}) \;\; \text{such that for all $1 \le j \le i$,} \\
	&\;\qquad X_{s_{j}} \in \oB(r_F) \;\; \text{and} \;\; \sU\Big(F_0,\sum_{m=1}^{j}s_m\Big) - \sU\Big(F_0, \sum_{m=1}^{j-1} s_m\Big) \ge  c_3 \lVert F \rVert_{1} \Theta(t_{n+1}/u_{n+1}) \bigg\}.
\end{align*}
Then by the strong Markov property, we get that for all $n \ge N$,
\begin{align}\label{e:occ_LIL_lower}
	\P^{o}\big(E_{n+1} \; | \; \sF_{t_n} \big)& \ge \P^{o}\big(E_{n+1} \cap G_n  \; | \; \sF_{t_n} \big) \nn\\[5pt]
	&\ge \inf_{w \in \oB(\vp^{-1}(n t_n))} \P^w \big(\,\sU(F, t_{n+1}-t_n) \ge \delta c_3 u_{n+1} \lVert F \rVert_{1} \Theta(t_{n+1}/u_{n+1}) \big) \cdot \1_{G_n}\nn\\
	&\ge \inf_{w \in \oB(\vp^{-1}(n t_n))} \P^w \big(\,\sU(F_0, t_{n+1}-t_n) \ge \delta c_3 u_{n+1} \lVert F \rVert_{1} \Theta(t_{n+1}/u_{n+1}) \big) \cdot \1_{G_n}\nn\\
	& \ge \inf_{w \in \oB(\vp^{-1}(nt_n))} \P^w \big(H_{n,\lfloor\delta u_{n+1}\rfloor +1}\big) \cdot \1_{G_n}.
\end{align}
Using the strong Markov property, \eqref{e:usize} and \eqref{e:occlower}, we get that for all $n$ large enough,
\begin{align}\label{e:occ_lower_onestep}
		\inf_{w \in \oB(\vp^{-1}(n t_n))} \P^w \big(H_{n,\lfloor\delta u_{n+1}\rfloor +1}\big)& = \inf_{w \in \oB(\vp^{-1}(n t_n))} \E \left[ \P^w \Big(H_{n,\lfloor\delta u_{n+1}\rfloor + 1 } \, \big| \, \sF_{\sum_{j=1}^{\lfloor\delta u_{n+1}\rfloor}s_j}\Big)\right]\nn\\[3pt]
	& \ge \inf_{w \in \oB(\vp^{-1}(n t_n))} \P^w \big(H_{n,\lfloor\delta u_{n+1}\rfloor }\big) \inf_{w \in \oB(r_F)} \P^w \big(H_{n,1 }\big) \nn\\
	&\ge  \inf_{w \in \oB(\vp^{-1}(n t_n))} \P^w \big(H_{n,\lfloor\delta u_{n+1}\rfloor }\big) \inf_{w \in \oB(\vp^{-1}(n t_n))} \P^w \big(H_{n,1 }\big)\nn\\
	& \ge\dots 
	\ge \Big(\inf_{w \in  \oB(\vp^{-1}(nt_n))} \P^w \big(H_{n,1}\big)\Big)^{\lfloor\delta u_{n+1}\rfloor+1 }.
\end{align}
Observe that for any $0<u<t$, if $\sU(F_0,t)-\sU(F_0,u)>0$, then there exists $s \in [u,t]$ such that $X_{s} \in {\rm supp}[F_0]\subset \oB(r_F)$. Thus,  $\inf_{w \in \oB(\vp^{-1}(nt_n))} \P^w(H_{n,1}) \ge p$ by  \eqref{e:occlower}. It follows from \eqref{e:occ_lower_onestep} and \eqref{e:usize} that for all $n$ large enough,
\begin{align*}
	\inf_{w \in \oB(\vp^{-1}(nt_n))} \P^w (H_{n,\lfloor\delta u_{n+1}\rfloor +1}) \ge p^{\delta u_{n+1} +1} \ge p^{3\delta \log(n+1)+1} = p(n+1)^{-1}.
\end{align*}
Hence, $\sum_{n=N}^\infty\inf_{w \in \oB(\vp^{-1}(nt_n))} \P^w (H_{n,\lfloor\delta u_{n+1}\rfloor +1}) = \infty$.

 Further, recall that  $t_{n+1} \ge 100t_n$ for all $n \ge N$. This implies that $t_n \ge 99^n$  for all $n$ large enough.
 Then by \eqref{e:scaling-infty}, it holds that for all $n$ large enough,
\begin{align*}
	\frac{2\vp^{-1}(nt_n)^{\up^{1/2}}}{\vp^{-1}(t_n)}  &= 	2(2R_\infty)^{\up^{1/2}-1} \bigg(\frac{\vp^{-1}(nt_n)}{\vp^{-1}(t_n)}\bigg)^{\up^{1/2}} 	\bigg(\frac{2R_\infty}{\vp^{-1}(t_n)}\bigg)^{1-\up^{1/2}} \\
	&\le c_4  n^{\up^{1/2}/\beta_1} t_n^{-(1-\up^{1/2})/\beta_2} \le c_4  n^{\up^{1/2}/\beta_1} 99^{-n(1-\up^{1/2})/\beta_2} <1.
\end{align*}
Therefore, we can apply Proposition \ref{p:EP}(ii) with $\up_1=\up^{1/2}$,  $r=\vp^{-1}(nt_n)$ and $t=t_n$  for all $n$ large enough and get that
\begin{align*}
\lim_{n \to \infty}	\P^o(G_n^c)\le \lim_{n \to \infty}  \P^o(\tau_{B(o, \vp^{-1}(nt_n))} \le t_n) \le \lim_{n \to \infty} c_5 n^{-1}=0.
\end{align*}

Finally, in view of \eqref{e:occ_LIL_lower}, we deduce from Lemma \ref{l:law001}  that $\P^z(\limsup E_n) =1$ and hence the lower bound in \eqref{e:occ_generalform_1} holds with $a_1=\delta c_3$.  The proof is complete. \qed	
\smallskip

\noindent \textbf{Proof of Corollary \ref{c:occ1}.} By \VDi \ and \eqref{e:scaling-infty*}, we see that for all $t>\vp(2R_\infty)$,
$$
\frac{t}{V(o, \vp^{-1}(t))}  = \frac{t}{ V(o,2R_\infty)} \frac{V(o,2R_\infty)}{V(o, \vp^{-1}(t))} \le c_1 t	\bigg( \frac{2R_\infty}{\vp^{-1}(t)} \bigg)^{d_1}  \le c_2t	\bigg( \frac{\vp(2R_\infty)}{t} \bigg)^{d_1/\beta_2}
$$
and
$$
	\frac{t}{V(o, \vp^{-1}(t))} \ge c_3 t	\bigg( \frac{2R_\infty}{\vp^{-1}(t)} \bigg)^{d_2}  \ge c_4t	\bigg( \frac{\vp(2R_\infty)}{t} \bigg)^{d_2/\beta_1}.
$$
Hence, since $\beta_1>d_2$, by Lemma \ref{l:theta}(iii), we see that \eqref{e:recur} holds with $x_0=o$ and
$$
\Theta(t)  \simeq\Theta(o, 2R_\infty, t) \simeq \frac{t}{V(o, \vp^{-1}(t))} \quad \text{and} \quad \log\log \Theta(t) \simeq  \log \log t \quad \;\;\text{for} \;\; t> 2 \vp(2R_\infty).
$$
Then using \VDi \ and \eqref{e:scaling-infty}, we deduce from  \eqref{e:occ_generalform_1} that there are constants $a_4 \ge a_3>0$ independent of $F$ such that
\begin{equation*}
	\limsup_{t \to \infty} \frac{\sU(F,t)/ \lVert F \rVert_{1}}{t/V\big(o, \vp^{-1}(t/\log \log t)\big)} \in [a_3,a_4],\quad\,
	\P^o\mbox{-a.s.}
\end{equation*}
For each $a>0$, one can see that $\left\{  \limsup_{t \to \infty} \frac{\sU(F,t)/\lVert F\rVert_1}{t/V(o, \vp^{-1}(t/\log \log t))} <a \right \}$ is shift-invariant by a similar proof to that  of Lemma \ref{l:occ1-1}.  Applying  Proposition \ref{p:law01} again, we conclude \eqref{e:occ_simpleform}. \qed

\section{Examples}\label{s:Example}
 The results of this paper cover a large class of subordinate diffusions and symmetric jump processes considered in  the authors' previous paper \cite{CKL}. Precisely, one can apply the  small time LIL result Theorem \ref{t:generalocc-SP}  to
 \cite[Examples 2.4, 2.5, 2.6(i), 2.8, 2.9]{CKL}, and  large time LIL results Theorem \ref{t:generalocc2} and one of Theorem \ref{t:generalocc3} or \ref{t:occ1} to  \cite[Examples 2.6(ii)]{CKL}.
  In this section, we give two important examples: (1) Feller processes on $\R^d$, and (2) Random conductance model. We begin with a lemma about stability of scale functions $\phi(x,r)$ and $\vp(r)$ in \eqref{e:occ_SP} and \eqref{e:generalocc2_1}. 
\begin{lem}\label{l:stability}
	Let $g$ and $h$ be increasing positive continuous functions  on a subinterval of $(0,\infty)$, and $\eps\in(0,1)$ be a constant.
	
	\noindent (i) If $\lim_{r \to 0}g(r)=0$, $g(r) \simeq h(r)$ for  $r\in (0,\eps)$ and there exist $x \in M$ and a constant $a_1\in (0,1]$ such that 
	$$\limsup_{r\to 0} \, \sU(B(x,r),g(r))/g(r)=a_1, \quad \P^x\text{-a.s.},$$
	then there exists a constant $a_2\in (0,1]$ such that 
		$$\limsup_{r\to 0} \, \sU(B(x,r),h(r))/h(r)=a_2, \quad \P^x\text{-a.s.}$$
		
	\noindent(ii) Suppose that zero-one law for shift-invariant event holds (see the paragraph above Proposition \ref{p:law01} for the definition of this) and $\lim_{r\to \infty}h(r+s)/h(r)=1$ for every $s>0$. If $\lim_{r \to \infty} g(r)=\infty$, $g(r) \simeq h(r)$ for $r>1/\eps$,  and there exist $x \in M$ and a constant $a_3\in (0,1]$ such that 
	\begin{equation}\label{e:stability-inf}
		\limsup_{r\to \infty} \, \sU(B(x,r),g(r))/g(r)=a_3, \quad \P^x\text{-a.s.},
	\end{equation}
	then there exists a constant $a_4\in (0,1]$ such that 
	$$\limsup_{r\to \infty} \, \sU(B(x,r),h(r))/h(r)=a_4,  \quad \P^z\text{-a.s.}, \;\; \forall z \in M.$$
\end{lem}
\pf Since the proofs are similar, we only give the proof for (ii).

As in the proof of Theorem \ref{t:generalocc2}, since zero-one law for shift-invariant event holds and  $\lim_{r\to \infty}h(r+s)/h(r)=1$ for every $s>0$, there exists a constant $a \in [0,1]$ such that 
	$$\limsup_{r\to \infty} \, \sU(B(x,r),h(r))/h(r)=a,  \quad \P^z\text{-a.s.}, \;\; \forall z \in M.$$
To conclude the result, it remains to show that $a>0$. Since $g(r) \simeq h(r)$ for $r>1/\eps$, there exists  $c_1\in (0,1]$ such that $c_1g(r)\le h(r) \le c_1^{-1}g(r)$ for $r>1/\eps$. Let $f(r):=g^{-1} (c_1g(r))$ which is well-defined for all $r$ large enough. Then $f(r) \le r$ since $g$ is increasing. Using \eqref{e:stability-inf} and the monotone property of occupation times, we get that   $\P^x$-a.s,
\begin{align*}
	a &= \limsup_{r\to \infty} \frac{\sU(B(x,r),h(r))}{h(r)}\ge  \limsup_{r\to \infty} \frac{\sU(B(x,r),c_1g(r))}{c_1^{-1}g(r)} = c_1^2  \limsup_{r\to \infty} \frac{\sU(B(x,r),g(f(r)))}{g(f(r))}\\
	& \ge c_1^2  \limsup_{r\to \infty} \frac{\sU(B(x,f(r)),g(f(r)))}{g(f(r))} = c_1^2  a_3>0.
\end{align*} 
This completes the proof. \qed

\subsection{LILs for Feller processes on $\R^d$}

Let $X$ be a (rich)  Feller process on $\R^d$ with the generator $(\sL, \sD(\sL))$ such that $C_c^\infty(\R^d) \subset \sD(\sL)$. It is well-known that $\sL$ restricted to $C_c^\infty(\R^d)$ is a \textit{pseudo-differential operator}, which has the following representation: For every $u \in C_c^\infty(\R^d)$,
\begin{align*}
	\sL u(x)&= -q(x,D)u(x):= -(2\pi)^{-d} \int_{\R^d} e^{i \la x,  \xi \ra}q(x, \xi) \int_{\R^d} e^{- i \la y, \xi\ra}u(y)dyd \xi
\end{align*}
with
\begin{align*}
	q(x, \xi)& = k(x) - i \la \b(x), \xi \ra +  \la \xi, \a(x)\xi \ra + \int_{\R^d \setminus \{0\}} (1- e^{i\la z, \xi \ra} + i \la z, \xi \ra \1_{\{| z| \le 1\}})\nu(x, dz)
\end{align*}
where $k:\R^d \to \R$ is  a non-negative measurable function,   $\b:\R^d \to \R^d$ a measurable function,  $\a:\R^d \to \R^{d \times d}$ a non-negative definite matrix-valued function, and $\nu(x,dz)$  a non-negative  $\sigma$-finite kernel on $\R^d \times \sB(\R^d \setminus \{0\})$ satisfying $\int_{\R^d \setminus \{0\}}(1 \land |z|^2) \nu(x, dz)<\infty$ for all $x \in \R^d$. The function $q:\R^d \times \R^d \to \mathbb C$ is called \textit{the symbol} of $X$ and the quadruplet $(k(x), \b(x), \a(x), \nu(x,dz))_{x \in \R^d}$ is called \textit{the L\'evy characteristics} of $X$.  We refer to \cite{Co65}. In this subsection, we always assume that $k(x)$ is identically zero so that $X$ has no killing inside.

Define for $x \in \R^d$ and $r>0$,
\begin{equation*}
	\Phi(x,r) = \bigg( \sup_{|\xi| \le 1/r} {\rm Re} \,\,q(x,\xi) \bigg)^{-1}.
\end{equation*}
We consider the following conditions on $X$. Fix an open set $U \subset \R^d$.
\smallskip

\noindent{\bf Assumption C.}
 Suppose that there exist constants $R_0, a_0, a_1,a_2,a_3,\eps_0\in (0,1)$ such that for any $x \in U$ and any $\xi \in \R^d$ with $|\xi|> 1/(R_0 \wedge (a_0 \updelta_U(x))$, 
 \begin{equation}\label{C1}\tag{C1}
\lim_{r \to 0} \Phi(x,r)=0 \quad  \text{(or, equivalently, either ${\bf a}(x) \neq 0$ or $\nu(x,\R^d)= \infty$)},
 \end{equation}\\[-7mm]
 \begin{equation}\label{C2}\tag{C2}
  \sup_{|\xi'|\le |\xi|}\text{Re} \, q(x, \xi') \ge a_1 |  \, \text{Im} \, q(x, \xi) |,
 \end{equation}\\[-6mm]
 \begin{equation}\label{C3}\tag{C3}
 \inf_{| x-y | \le 1/|\xi|} \textrm{Re}\,q(y,\xi) \ge a_2 \sup_{| x-y|\le 1/|\xi|} \textrm{Re}\,q(y,\xi),
 \end{equation}\\[-6mm]
 \begin{equation}\label{C4}\tag{C4}
  \inf_{z \in \R^d, \, |z| =1 }\P^x\big(2\la  X_t- x, z \ra \le -| X_t -x |  \big) \ge a_3 \quad \text{for  all  $0<t < \eps_0\Phi\big(x, R_0 \wedge (a_0 \updelta_U(x))\big)$}.
 \end{equation}
 
\begin{remark}
	{\rm When $d=1$ and $X$ is symmetric, \eqref{C4} holds true with $a_3=1/2$. }
\end{remark}

 From Theorem \ref{t:generalocc-SP}, we obtain the following limsup LIL for occupation times of $X$.  

\begin{prop}\label{p:Feller2}
	Let $X$ be a Feller process on $\R^d$ with symbol $q$. Suppose that Assumption C holds for an open set $U \subset \R^d$. Let $f_0:U \times (0,\infty) \to [0,1]$ be a deterministic function such that for every $x \in U$ and $\kappa>0$,
	\begin{equation}\label{e:ex-Feller} \limsup_{r \to 0} \frac{ \sU\big(B(x,r), \kappa \Phi(x,r) \log |\log \Phi(x,r)|\big)}{\kappa \Phi(x,r) \log |\log \Phi(x,r)|}=f_0(x,\kappa), \quad \P^x\mbox{-a.s.} 
	\end{equation}
	Then  $f_0$ satisfies all properties {\rm (P1)-(P3)} in Theorem \ref{t:generalocc-SP}.
\end{prop}
\pf Following the arguments in the proof of \cite[Theorem 2.4]{CKL2}, we deduce that if Assumption C holds true for $U$ then so \ref{A1}, \ref{A2}, \ref{A3} and \ref{A4} do for $U$. Now the result follows from Theorem \ref{t:generalocc-SP}.   \qed

A L\'evy process on $\R^d$  is a Feller process whose characteristics is independent of $x \in \R^d$. Hence,  $q(x,\xi)=\psi(\xi)$ for a negative definite function $\psi$ when $X$ is a L\'evy process. The function $\psi$ is called \textit{the L\'evy exponent}. For a given L\'evy exponent $\psi$, we write
\begin{equation}\label{d:scale-Levy}
 \Phi_\psi(r):=\Big( \sup_{|\xi| \le 1/r} {\rm Re} \,\,\psi(\xi) \Big)^{-1}.
\end{equation}

\begin{example}\label{E:Isotropic Levy}
	{\rm {\bf (Isotropic L\'evy processes)} Let $X$ be a  L\'evy process on $\R^d$ with exponent  $\psi$. Suppose  that  $\psi$ is a radial function, namely, $\psi(\xi)=\psi^*(|\xi|)$ for some $\psi^*$ and $\lim_{r \to \infty} \psi^*(r)=\infty$. Here we do not assume the absolute continuity of the L\'evy measure $\nu$.
		
	Note that 	 $\overline{\psi(\xi)} = \psi(-\xi) = \psi(\xi) = \text{Re} \, \psi(\xi)$ for all $\xi \in \R^d$.  Hence, \eqref{C1} and \eqref{C2} are satisfied for $U=\R^d$. Clearly, \eqref{C3} holds true for $U=\R^d$. Lastly, we get \eqref{C4} for $U=\R^d$ from \cite[Proposition 5.2]{GRT} and \cite[Lemma 2.2]{CKL2}. Therefore, we can apply Proposition  \ref{p:Feller2} and deduce that there exist constants $c_1,\kappa_1>0$ and a non-increasing function function $f_0$ on $(0,\infty)$ satisfying $f_0(\kappa)=1$ for $\kappa \le \kappa_1$   and  $\lim_{\kappa \to \infty} f_0(\kappa)=0$ such that  for every $x \in \R^d$ and $\kappa>0$,
		\begin{equation}\label{ex:Levy2} \limsup_{r \to 0} \frac{ \sU\big(B(x,r), \kappa \Phi_\psi(r) \log |\log \Phi_\psi(r)|\big)}{\kappa \Phi_\psi(r) \log |\log \Phi_\psi(r)|}=f_0(\kappa), ~\quad \P^x\mbox{-a.s.},
		\end{equation}
		where the function $\Phi_\psi$ is defined as \eqref{d:scale-Levy}. 
		
		Note that the LIL \eqref{ex:Levy2} covers L\'evy processes  whose exponent is slowly varying at infinity.
		  \qed
	}
\end{example}

In \cite{CKL2}, the authors give a sufficient condition for Assumption C in terms of the symbol $q(\cdot,\xi)$ only by using the symmetrization argument from  \cite{SW13}. Fix an open set $U \subset \R^d$ as before.

\medskip

\noindent {\bf Assumption S.}  $C_c^\infty(\R^d)$ is an operator core for  $(\sL, \sD(\sL))$, i.e.   $\overline{\sL|_{C_c^\infty(\R^d)}}=\sL$, and there exist  constants  $R_0,  a_0\in (0,1)$, $c_L,c_U>0$ and $K\ge 1$ such  that the following properties are satisfied  for every $x \in U$:

\smallskip

\begin{enumerate}[(S1)]
	\item There exists an increasing function $g(x,\cdot)$ and constants $0<\alpha(x)\le \beta(x)$  such that 
	\begin{equation}\label{a:alphabeta}
		\frac{1}{\alpha(x)}-\frac{1}{\beta(x)} < \frac{1}{d^2+d},
	\end{equation}
	$$
c_L \bigg(\frac{r}{s}\bigg)^{\alpha(x)}\le \frac{g(x,r)}{g(x,s)} \le c_U \bigg(\frac{r}{s}\bigg)^{\beta(x)} \quad \text{for all} \;\, r \ge s> 1/(R_0 \wedge (a_0 \updelta_U(x)))
$$
and
	\begin{equation}\label{a:as1}
	\hspace{15mm}	K^{-1}g(x,|\xi|) \le \text{Re}\, q(x, \xi) \le K g(x,|\xi|)\; \text{ for all } \, \xi \in \R^d, \; |\xi|>1/(R_0 \wedge (C_0 \updelta_U(x))).
	\end{equation}

	\item For every $0<R<R_0 \wedge (C_0 \updelta_U(x))$, there exists  a Feller process $Y=Y^{x,R}$ with symbol $q_Y(\cdot,\xi)$ such that
	\begin{align}
		&(i) \;  q(y, \xi)=2 \text{Re}\,q_Y(y,\xi/2) \; \text{ for all } y \in \overline{B(x,R)} \text{ and } \xi\in \R^d,\label{a:as2}\\[5pt]
		&(ii) \; K^{-1} \inf_{|z-x| \le r} \textrm{Re}\,q_Y(z,\xi) \le  \textrm{Re}\,q_Y(y,\xi) \le K \sup_{|z-x|\le r} \textrm{Re }q_Y(z,\xi) \nn\\
		&\quad\;\; \text{for all } y \in \R^d \setminus \overline{B(x,R)} \text{ and } \xi\in \R^d, \; |\xi|> 1/(R_0 \wedge (C_0 \updelta_U(x))). \label{a:as3}
	\end{align}
\end{enumerate}

\bigskip

We refer to \cite[Examples 2.13 and 2.14]{CKL2} for concrete examples of Feller processes satisfying Assumption S.

 \begin{remark}\label{r:well-posed}
 	{\rm $C_c^\infty(\R^d)$ is an operator core for $(\sL, \sD(\sL))$ if and only if the martingale problem for $(-q(\cdot, D), C_c^\infty(\R^d))$ is well-posed.	See \cite[Proposition 4.6]{SW13}.
 	}
 \end{remark} 
 
 From Proposition \ref{p:Feller2} and \cite[Proposition 2.12]{CKL2}, we conclude 
\begin{cor}
 Let $X$ be a Feller process on $\R^d$ with symbol $q$. Suppose that \eqref{C3} and Assumption S hold true for an open set $U \subset \R^d$. Then the deterministic function $f_0:U \times (0,\infty) \to [0,1]$ defined by \eqref{e:ex-Feller}   satisfies all properties {\rm (P1)-(P3)} in Theorem \ref{t:generalocc-SP}.
\end{cor}

\subsection{Random conductance model}\label{s:RCM}

 In this section, we give LILs for occupations times  at infinity for some random conductance models.  We repeat the setting of the  random conductance models in  \cite[Section 3]{CKL} here for the readers' convenience. 

Let $G= (\LL, E_\LL)$ be a locally finite connected infinite undirected graph, where $\LL$ is the set of vertices and $E_\LL$ the set of edges. We denote $d(x,y)$ for the graph distance between $x,y \in \LL$, and  $\mu_c$ for the counting measure on $\mathbb{L}$.  A \textit{random conductance} $\boldsymbol{\eta}=(\eta_{xy} : x,y \in \mathbb{L})$ on $\mathbb{L}$ is a family of nonnegative random variables defined on some probability space $(\boldsymbol\Omega,{\bf F}, {\bf P})$ such that  $\eta_{xx} = 0$ and $\eta_{xy} = \eta_{yx}$ for all $x,y \in \mathbb{L}$. 
Associated with a random conductance $\boldsymbol{\eta}$, for each $\boldsymbol \omega \in \boldsymbol\Omega$, the \textit{variable speed random walk} (VSRW) $X^{{\boldsymbol \omega}} = (X^{\boldsymbol \omega}_t, t \ge 0; {\mathbb P}^x_{\boldsymbol \omega}, x \in \mathbb{L})$ is defined by the symmetric Markov process on $\mathbb{L}$ with $L^2(\mathbb{L}, \mu_c)$-generator
$$ \sL_V^{\boldsymbol \omega} f(x) = \sum_{y \in \mathbb{L}} \eta_{xy}({\boldsymbol \omega})(f(y) - f(x)),  $$
and the \textit{constant speed random walk} (CSRW) $Y^{\boldsymbol \omega} = (Y^{\boldsymbol \omega}_t, t \ge 0; {\mathbb P}^x_{\boldsymbol \omega}, x \in \mathbb{L})$  is the symmetric Markov process on $\mathbb{L}$ with $L^2(\mathbb{L}, \nu)$-generator
$$ \sL_C^{\boldsymbol \omega} f(x) = \nu_x^{-1}({\boldsymbol \omega}) \sum_{y \in \mathbb{L}} \eta_{xy}({\boldsymbol \omega})(f(y)-f(x)), $$
where $\nu_x:=\sum_{y \in \mathbb{L}} \eta_{xy}$.  

In this section, we denote by $\sU_X(F,t)$ and $\sU_Y(F,t)$ the occupation times of VSRW $X$ and CSRW $Y$ for $F$ respectively.

\subsubsection{Random walks on supercritical bond percolation clusters}

\

\smallskip

Let $\LL=\Z^d$, $d \ge 2$ and $p \in (0,1]$. Let $(\eta_{xy})_{x,y \in \Z^d}$ be independent random variables  such that
\begin{equation}\label{e:percol}
	\eta_{xy}= \begin{cases}
		\text{Bernoulli distribution with mean } p, & \text{if} \;\; |x-y|=1, \\
		0, & \text{otherwise},
	\end{cases}
\end{equation}
An edge $\{x,y\}$ in $\Z^d$ is called \textit{open} if $\eta_{xy}=1$ and a set $\sC \subset \Z^d$ is called an \textit{open cluster} if every $x, y \in \sC$ are connected by an open path.  It is known that there exists a critical value $p_c=p_c(d)\in (0,1)$ such that if $p>p_c$, then $\bf P$-a.s., there exists  a unique infinite open cluster in each configuration, which we denote $\sC_\infty = \sC_\infty({\boldsymbol \omega})$. It is obtained in \cite{Keper} that $p_c(2)=1/2$, but the exact values of $p_c(d)$ for $d\ge 3$  are still open. For details and a comprehensive bibliography,  we refer to   \cite{Ke2, Gri}.

Recall the definition of  $\fMg$ from \eqref{e:Linfty}. Below, we give LILs for  occupation times of both VSRW $X$ and CSRW $Y$ on $\sC_\infty$ associated with  the random conductance $\boldsymbol \eta$ defined by \eqref{e:percol}.

\begin{prop}\label{p:rcm-1}
	Let $p \in (p_c,1]$. There exist $\boldsymbol \Omega_0\subset \boldsymbol \Omega$ with ${\bf P}({\boldsymbol \Omega_0})=1$ and constants $\kappa_3 \ge \kappa_2 \ge \kappa_1 >0$, $a_2 \ge a_1>0$ which depend only on the parameter $p$  such that  the following statements hold true  for every ${\boldsymbol \omega} \in \boldsymbol \Omega_0$.
	
	\noindent (i) There exists an non-increasing function $f_\infty(\cdot, {\boldsymbol \omega}):(0,\infty) \to (0,1]$ such that $f_\infty(\kappa,{\boldsymbol \omega})=1$ if $\kappa \le \kappa_1$, $\lim_{\kappa \to \infty} f_\infty(\kappa, {\boldsymbol \omega})=0$ uniformly on ${\boldsymbol \omega} \in \boldsymbol \Omega_0$ and that for every $\kappa>0$, 
	\begin{equation*}
		\limsup_{r \to \infty} \frac{\sU_X\big(B(x, r),\, \kappa r^2 \log \log r\big)}{\kappa r^2 \log \log r} = f_\infty(\kappa,{\boldsymbol \omega}),\qquad \P^z_{{\boldsymbol \omega}}\text{-a.s.}, \;\; \forall z \in \sC_\infty({\boldsymbol \omega}).
	\end{equation*}
	
	\noindent (ii) (a) If $d \ge 3$, then there exists a constant $\kappa({\boldsymbol \omega})\in [\kappa_2,\kappa_3]$ such that 
	\begin{equation*}
	\limsup_{r \to \infty} \frac{\sU_X\big(B(x, r),\, \infty \big)}{r^2 \log \log r} = \kappa({\boldsymbol \omega}),\qquad \P^z_{{\boldsymbol \omega}}\text{-a.s.}, \;\; \forall z \in \sC_\infty({\boldsymbol \omega}).
\end{equation*}

	(b) If $d=2$, then  for every  $F \in \mathfrak B_{b,+}^\infty(\sC_\infty({\boldsymbol \omega});\mu_c) $ with $\lVert F \rVert_{L^1(\sC_\infty({\boldsymbol \omega}); \mu_c)} \neq 0$, there exists a constant $a_F({\boldsymbol \omega}) \in [a_1, a_2]$ such that 
	\begin{equation}\label{e:2.14-example}
		\limsup_{t \to \infty} \frac{\sU_X(F, t)/\lVert F \rVert_{L^1(\sC_\infty({\boldsymbol \omega});\mu_c)}}{\log t\, \log \log \log t} = a_F, \qquad \P^z_{{\boldsymbol \omega}}\text{-a.s.}, \;\; \forall z \in \sC_\infty({\boldsymbol \omega}).
	\end{equation}

\noindent (iii) The assertions in (i-ii)  hold true for $\sU_Y$ instead of $\sU_X$ with different constants $\kappa_1, \kappa_2, \kappa_3,a_1,a_2$.
\end{prop}
\pf   According to \cite[Theorem 2.18(c), Lemma 2.19, Theorem 5.7 and Lemma 5.8]{Ba} and the Borel-Cantelli lemma,  for $\bf P$-a.s. ${\boldsymbol \omega}$,  there is a constant $R_\infty({\boldsymbol \omega})\in [1,\infty)$ such that  $\mathrm{VRD}^{R_\infty({\boldsymbol \omega})}(4/5)$ with respect to $\mu_c$ holds, and $\mathrm{NDL}^{R_\infty}(\varPhi,4/5)$ and $\mathrm{NDU}^{R_\infty}(\varPhi,4/5)$ are satisfied with $\varPhi(r)=r^2$.  $\mathrm{Tail}^{R_\infty}(\varPhi,4/5,\le)$ trivially holds because there are no long range jumps. Therefore, the assertions in (i) and (ii) are consequences of Theorems \ref{t:generalocc2}, \ref{t:generalocc3} and \ref{t:occ1}.  Since $\mu_c|_{\sC_\infty} \le \nu_\cdot |_{\sC_\infty} \le 4\mu_c|_{\sC_\infty}$, we also obtain (iii). \qed

\subsubsection{Random conductance model with unbounded conductances}

\

\smallskip

Let $\LL=\Z^d$, $d \ge 2$ and $\xi$ be any random variable such that $\xi \ge 1$, $\bf P$-a.s. Let $(\eta_{xy})_{x,y \in \Z^d}$ be independent random variables such that
\begin{align}\label{e:unbdd}
	\eta_{xy}= \begin{cases}
		\text{an independent copy of } \xi, & \text{if} \;\; |x-y|=1, \\
		0, & \text{otherwise}.
	\end{cases}
\end{align}	
We establish LILs for occupation times of VSRW $X$ associated with \eqref{e:unbdd}.

\smallskip

\begin{prop}
Assertions (i) and (ii) of Proposition \ref{p:rcm-1} hold true for $\sU_X(\cdot, \cdot)$ with $\Z^d$ instead of $\sC_\infty({\boldsymbol \omega})$, where $X$ is the VSRW associated with the random conductance defined by \eqref{e:unbdd}.
\end{prop}
\pf Clearly $\mathrm{VRD}^2(4/5)$ holds for $(\Z^d, \mu_c)$. Using \cite[Lemma 3.5 and Theorems 4.3 and 4.6]{BD} and the Borel-Cantelli lemma, we deduce that there exists for $\bf P$-a.s. ${\boldsymbol \omega}$, there is a  constant $R_\infty({\boldsymbol \omega})\in(1,\infty)$ such that  $\mathrm{NDL}^{R_\infty}(\varPhi,4/5)$ and $\mathrm{NDU}^{R_\infty}(\varPhi,4/5)$  hold with $\varPhi(r)=r^2$. Since  $\mathrm{Tail}^{R_\infty}(\varPhi,4/5,\le)$ clearly holds,  using Theorems \ref{t:generalocc2}, \ref{t:generalocc3} and \ref{t:occ1}, we get the result.\qed

\subsubsection{Random conductance model with stable-like jumps}

\

\smallskip

In  \cite{CKL, CKL2}, the authors have obtained limsup and liminf LILs at infinity for random conductance model with long range jumps using results from \cite{CKW18, CKW20-1}. Here we give LILs for occupation times at infinity for this model.

 Suppose that there is a constant $d>0$ such that
\begin{equation}\label{(3.1)}
	\mu_c(B(x,r)) \simeq r^d \quad \mbox{for all} \;\, x \in \mathbb{L}, \, r>10
\end{equation}
Let  $\alpha \in (0 \vee (2-\frac d2), 2)$ and  $\boldsymbol \eta$ be a random conductance on $\LL$  such that $w_{xy}:=\eta_{xy}|x-y|^{d+\alpha}$ satisfies the following properties:
\begin{equation}\label{e:long-range-jump}
	\sup_{x,y \in \LL, \,x \neq y} {\bf  P}(w_{xy}=0)<1/2 \qquad \text{and} \;\;\quad \sup_{x,y \in \LL, \,x \neq y} {\bf E} \big[w_{xy}^p + w_{xy}^{-q}\1_{\{w_{xy}>0\}}\big]<\infty, 
\end{equation}
where the constants $p$ and $q$ satisfy
$$
p>\frac{d+2}{d} \vee  \frac{d+1}{4-2\alpha}, \qquad q>\frac{d+2}{d}.
$$
When we consider the CSRW $Y^{\boldsymbol \omega}$, we additionally assume that there exist constants $m_2 \ge m_1>0$ such that for ${\bf P}$-a.s. ${\boldsymbol \omega}$,
\begin{equation}\label{e:CSRW}
\eta_{xy}({\boldsymbol \omega})>0 \text{ for all $x,y \in \LL, \, x \neq y$ \; and \;} m_1 \le \sum_{y \in \LL} \eta_{xy}({\boldsymbol \omega}) \le m_2 \text{ for all $x \in \LL$}.
\end{equation}

\begin{prop}\label{p:rcm-3}
 Suppose that \eqref{(3.1)} and \eqref{e:long-range-jump} are satisfied. There exist  ${\boldsymbol \Omega_0}\subset \boldsymbol \Omega$ with ${\bf P}({\boldsymbol \Omega_0})=1$ and constants $\kappa_3 \ge \kappa_2 \ge \kappa_1 >0$, $a_2 \ge a_1>0$ such that  the following statements hold true for every ${\boldsymbol \omega} \in {\boldsymbol \Omega_0}$.
	
	\noindent (i) There exists an non-increasing function $f_\infty(\cdot, {\boldsymbol \omega}):(0,\infty) \to (0,1]$ such that $f_\infty(\kappa,{\boldsymbol \omega})=1$ if $\kappa \le \kappa_1$, $\lim_{\kappa \to \infty} f_\infty(\kappa, {\boldsymbol \omega})=0$ uniformly on ${\boldsymbol \omega} \in {\boldsymbol \Omega_0}$ and that for every $\kappa>0$, 
	\begin{equation*}
		\limsup_{r \to \infty} \frac{\sU_X\big(B(x, r),\, \kappa r^\alpha \log \log r\big)}{\kappa r^\alpha \log \log r} = f_\infty(\kappa,{\boldsymbol \omega}),\qquad \P^z_{{\boldsymbol \omega}}\text{-a.s.}, \;\; \forall z \in \LL.
	\end{equation*}
	
	\noindent (ii)(a) If $d < \alpha$, then  for every $\gamma \in (\frac{\alpha d}{\alpha-d},\infty]$ and $F \in \mathfrak B_{b,+}^\gamma(\LL;\mu_c) $ with $\lVert F \rVert_{L^1(\LL; \mu_c)} \neq 0$, there exists a constant $a_F({\boldsymbol \omega}) \in [a_1, a_2]$ such that 
	\begin{equation*}
		\limsup_{t \to \infty} \frac{\sU_X(F, t)/\lVert F \rVert_{L^1(\LL;\mu_c)}}{t^{1-1/\alpha} (\log \log t)^{1/\alpha}} = a_F, \qquad \P^z_{{\boldsymbol \omega}}\text{-a.s.}, \;\; \forall z \in \LL.
	\end{equation*}
	
	(b) If $d >\alpha $, then there exists a constant $\kappa({\boldsymbol \omega})\in [\kappa_2,\kappa_3]$ such that 
	\begin{equation*}
		\limsup_{r \to \infty} \frac{\sU_X\big(B(x, r),\, \infty \big)}{\kappa({\boldsymbol \omega}) \,r^\alpha \log \log r} = 1,\qquad \P^z_{{\boldsymbol \omega}}\text{-a.s.}, \;\; \forall z \in \LL.
	\end{equation*}

\noindent (iii) If \eqref{e:CSRW} are also satisfied, then the assertions in (i-ii)  hold true for $\sU_Y$ instead of $\sU_X$ with different constants $\kappa_1, \kappa_2, \kappa_3,a_1,a_2$.
\end{prop}
\pf Note that \eqref{(3.1)} yields  $\mathrm{VRD}^{10}(\up)$ with respect to $\mu_c$ for every $\up \in (0,1)$ and \eqref{e:CSRW} yields $\mathrm{VRD}^{10}(\up)$ with respect to $\nu_\cdot$. By  the proof of \cite[Theorem 3.1]{CKL} and \cite[Proposition 2.3]{CKW20-1},  there exist $\up\in(0,1)$  such that for $\bf P$-a.s. ${\boldsymbol \omega}$, with some constant $R_\infty({\boldsymbol \omega})\in (1,\infty)$, \ref{B4+} and \ref{NDU} are satisfied with $\varPhi(r)=r^\alpha$.
 Applying  Theorems \ref{t:generalocc2}, \ref{t:generalocc3} and \ref{t:occ1}, we arrive at the desired result. \qed

\begin{remark}
{\rm The assertions of Proposition \ref{p:rcm-3}  hold true for $\alpha \in (0 \vee (1- \frac d2),1)$, if \eqref{e:long-range-jump} is satisfied with some $p> \frac{d+2}{d} \vee \frac{d+1}{2-2\alpha}$ and $q>\frac{d+2}{d}$, by the proof for the second statement of \cite[Theorem 3.1]{CKL}.
}
\end{remark}

\section{Appendix}

In this appendix, we give two lemmas which are used in the proof of Theorem \ref{t:occ1}.

The first lemma is a modification  of \cite[Theorem 3.1]{BK} which gives upper bounds for LILs. 

\begin{lem}\label{l:limsup}
	Let $\X$ be a topological space and $Y=(Y_t, t \ge 0;\, \P^y, y \in \X_\partial)$ a strong Markov process on $\X$ with a cemetery point $\partial$. Denote by $\theta^Y_s$ the shift operator with respect to $Y$.
	Let $(F_t)_{t \ge 0}$ be a continuous non-decreasing adapted functionals of $Y$ and $g:(0,\infty) \to (0,\infty)$ be  an increasing continuous function.
	Suppose that there exists  $T>0$ such that the following conditions are satisfied:

 {\rm (1)} $F(t) - F(s) \le F(t-s) \circ \theta^Y_s\,$ for all  $T < s \le t$.
	
	{\rm (2)} $\lim_{t \to \infty} g(t) = \infty$ and there exists a constant $a_1>1$ such that
	\begin{equation}\label{e:condition_g2}
		g(2t) \le a_1g(t) \quad \text{for all} \;\, t>T.
	\end{equation}

	{\rm (3)} There exist constants $a_2, \eps_1>0$ such that 
	\begin{equation}\label{e:condition_g1}
		\sup_{z \in \X} \P^z \big( F(t) \ge a_2g(t) \big) \le e^{-\eps_1} \quad \text{for all} \;\, t>T.
	\end{equation}

	\noindent Then, there exists a constant $a_3 \in (0,\infty)$ which depends only on $a_1,a_2$ and $\eps_1$ such that
	$$
	\limsup_{t \to \infty} \frac{ F(t) }{g\big(t/\log \log g(t) \big) \log \log g(t) } \le a_3, \quad\;\; \P^z\mbox{-a.s.},~~\forall z\in \X.
	$$
\end{lem}
\pf Let $K:=\log a_1$. Since $\lim_{t \to \infty}g(t) = \infty$, there exist $N>3$ and  a sequence $(t_n)_{n \ge N}$ such that 
\begin{equation}\label{e:deftn}
	g\big(t_n/ \log (nK)\big) = e^{nK} \qquad \text{for all} \;\; n \ge N.
\end{equation}
By \eqref{e:condition_g2},  it holds that
$$g\big(2t_{n}/\log (nK)\big) \le a_1 e^{nK} \le e^{(n+1)K} = g\big(t_{n+1}/\log (nK+K)\big).$$ Since $g$ is increasing, this yields that
\begin{equation}\label{e:ratio}
t_{n+1} \ge 2t_n \quad \text{for all } \;n \ge N.
\end{equation}
Define $u_n=\log \log g(t_n)$ for $n \ge N$. By \eqref{e:deftn} and \eqref{e:condition_g2}, we have 
\begin{align*}
	e^{nK} \le g(t_n) \le c_1(\log (nK))^{c_2} g(t_n/\log(nK)) = c_1(\log (nK))^{c_2}  e^{nK}.
\end{align*}
Since $t_n$ grows at least exponentially by \eqref{e:ratio}, we deduce that there exists $N' \ge N$ such that
\begin{equation}\label{loglogg2}
\log n \le u_n \le 2\log n <t_n/T \qquad  \text{for all } \; n \ge N'.
\end{equation}
Applying the property (1) yields that for every $n \ge N'$,
\begin{align}\label{decomp}
 F(t_{n}) \le\sum_{j=0}^{\lfloor u_n \rfloor}   F(t_n/u_n)\circ \theta^Y_{jt_n/u_n}.
\end{align}
We claim that for all $n \ge N'$ and $m \ge 1$,
\begin{align}\label{induction}
	\sup_{z \in \X}\P^z \big(  F(t_n/u_n)\ge a_2m g( t_n/ u_n) \big) \le e^{-m \eps_1}.
\end{align}
	By \eqref{e:condition_g1},  \eqref{induction} holds true when $m=1$.
Assume that  \eqref{induction} holds for   $m \le k\in \N$. Set 
$$T_1:=\inf \big\{s > 0: F(s) \ge a_2g(t_n/ u_n) \big\}.$$ 
Using the strong Markov property, \eqref{e:condition_g1} and the induction hypothesis, we get that for all $z \in \X$,
\begin{align*}
	\P^z \big(F(t_n/u_n) \ge a_2(k+1)g( t_n/ u_n) \big) &= \P^z  \big( F(t_n/u_n) \ge a_2(k+1)g(t_n/u_n), \; T_1 \le  t_n/u_n \big) \nn\\
	&\le \E^z \P^{X_{T_1}} \big( F(t_n/u_n)  \ge a_2 k g(t_n/u_n) \big) \P^z \big( T_1 \le  t_n/u_n\big) \le e^{-(k+1)\eps_1}.
\end{align*}
Hence,  we conclude that \eqref{induction} is true by induction.
 
 Using Markov inequality and \eqref{induction}, we arrive at
\begin{align}\label{e:exponent-bound}
	&	\sup_{z \in \X}\E^z \left[ \exp\bigg( \frac{\eps_1}{2a_2}\frac{ F(t_n/u_n) }{g(t_n/ u_n)} \bigg) \right] \nn\\
	&\le \sum_{m=0}^\infty e^{\eps_1 (m+1)/2} 	\sup_{z \in \X}\P^z\left( \frac{F(t_n/u_n)}{g( t_n/ u_n)}  \in \big[a_2m, a_2(m+1)\big)\right) \le e^{\eps_1/2}\sum_{m=0}^\infty e^{-m \eps_1/2} = c_3.
\end{align}
By Markov inequality, \eqref{loglogg2} and \eqref{decomp}, \eqref{e:exponent-bound} yields that  for all $z \in \X$ and $n \ge N'$,
\begin{align*}
	&\P^z\big( F(t_n) > 2a_2\eps_1^{-1}(3+\log c_3) g(t_n/u_n) u_n\big) \le e^{-(3+\log c_3) u_n} \E^z \left[ \exp\bigg( \frac{\eps_1}{2a_2}\frac{F(t_n)}{g(t_n/u_n)} \bigg) \right] \nn\\
	& \le e^{-(3+\log c_3) u_n}\prod_{j=0}^{\lfloor u_n \rfloor} \sup_{w \in \X} \E^w \left[ \exp\bigg(\frac{\eps_1}{2a_3}  \frac{F(t_n/u_n ) }{g(t_n/u_n)} \bigg) \right]  \le c_3 e^{-3u_n}  \le c_3n^{-3/2}.
\end{align*}
Finally, using the Borel-Cantelli lemma, \eqref{e:condition_g2}, \eqref{e:deftn} and \eqref{loglogg2}, it holds that for all $z \in \X$ and $\P^z$-a.s.,
\begin{align*}
	&\limsup_{t \to \infty} \frac{ F(t) }{g\big(t/\log \log g(t)\big) \log \log g(t) }\\
	&\le \limsup_{n \to \infty} \sup_{t \in [t_{n-1},t_{n}]} \frac{ F(t) }{g\big(t/\log \log g(t)\big) \log \log g(t) } \le \limsup_{n \to \infty} \frac{ F(t_n) }{g\big(t_{n-1}/u_{n-1}\big) u_{n-1} }\\
	&\le c_4\limsup_{n \to \infty} \frac{ F(t_n) }{u_{n-1}e^{(n-1)K}  } \le c_5e^K \limsup_{n \to \infty} \frac{ F(t_n) }{g\big(t_{n}/u_{n}\big) u_{n} } \le  2a_2\eps_1^{-1}(3+\log c_3) c_5e^K.
\end{align*}
The proof is complete. \qed

The second lemma is a version of conditional Borel-Cantelli lemma which was proved in \cite[Proposition 5.1]{CKL}.

\begin{lem}\label{l:law001}
	Let $(\tO, \tP, \GG, (\GG_t)_{t \ge 0})$ be a filtered probability space and $(s_n)_{n\ge 1}$ an increasing sequence with $\lim_{n \to \infty} s_n = \infty$.  Assume that families of events $(E_n)_{n \ge 1}$ and $(G_n)_{n \ge 1}$ satisfy the following properties:

	{\rm (1)} $E_n, G_{n} \in \GG_{s_n}$ for all $n \ge 1$.
	
 {\rm (2)} There exist a constant $p>0$ and  a sequence $(b_n)_{n \ge 1}$ with $\sum_{n=1}^\infty b_n = \infty$ such that
	\begin{align*}
		\tP(G_n^c) \le p \quad \text{and} \quad \tP(E_{n+1} \, | \, \GG_{s_{n}}) \ge b_{n} \1_{G_{n}} \quad \text{for all $n \ge 1$.}
	\end{align*}
	
	\noindent	Then,  $\tP(\limsup E_n) \ge 1-p$. In particular, if $\lim_{n \to \infty} \tP(G_n^c)=0$, then $\tP(\limsup E_n) = 1$.	
\end{lem}

\end{document}